\documentclass[11pt,twoside]{amsart}
\usepackage{amscd,moduli}

\author{Valery Alexeev and Michel Brion}
\address{Department of Mathematics\\
University of Georgia\\
Athens, GA 30602, USA}
\email{valery@math.uga.edu}
\address{Institut Fourier, B. P. 74\\
38402 Saint-Martin d'H\`eres Cedex, France}
\email{Michel.Brion@ujf-grenoble.fr}

\begin{document}

\begin{abstract}
For a connected reductive group $G$ and a finite-dimensional
$G$-module $V$, we study the invariant Hilbert scheme that
parameterizes closed $G$-stable subschemes of $V$ affording a 
fixed, multiplicity-finite representation of $G$ in their 
coordinate ring. We construct an action on this invariant 
Hilbert scheme of a maximal torus $T$ of $G$, together with 
an open $T$-stable subscheme admitting a good quotient. 
The fibers of the quotient map 
classify affine $G$-schemes having a prescribed categorical 
quotient by a maximal unipotent subgroup of $G$. We show that $V$ 
contains only finitely many multiplicity-free $G$-subvarieties, 
up to the action of the centralizer of $G$ in $\GL(V)$. As a 
consequence, there are only finitely many isomorphism classes 
of affine $G$-varieties affording a prescribed multiplicity-free 
representation in their coordinate ring.
\end{abstract}

\bibliographystyle{amsalpha}

\title{Moduli of affine schemes with reductive group action}

\date{June 17, 2003}

\maketitle 

\tableofcontents

\setcounter{section}{-1}

\section{Introduction}

A fundamental result in the classification theory of projective 
algebraic varieties is the existence of the Hilbert scheme,
a projective scheme parameterizing closed subschemes of a
fixed projective space, having a fixed Hilbert polynomial.
In this paper, we study two versions of the Hilbert scheme,
relevant to the classification of affine algebraic varieties
$X$ equipped with an action of a reductive group $G$. 

In this setting, an analog of the Hilbert polynomial is 
the representation of $G$ in the coordinate ring of $X$.
This rational, possibly infinite-dimensional $G$-module 
is the direct sum of simple modules, with multiplicities 
that we assume to be all finite. Now the closed 
$G$-stable subschemes of a fixed finite-dimensional 
$G$-module $V$, having fixed multiplicities in their 
coordinate ring, are parameterized by a quasi-projective 
scheme: the invariant Hilbert scheme $\Hilb^G_h(V)$ 
(where $h$ encodes the multiplicities).

For subschemes of finite length, the invariant Hilbert scheme 
may be realized in a punctual Hilbert scheme $\Hilb_n(V)$, 
as a union of connected components of the $G$-invariant 
subscheme. This construction is well known in the case of 
finite groups, see e.g. \cite{Nakamura01}. But its generalization 
to arbitrary subschemes is problematic, since there is no 
reasonable moduli space for closed subschemes of infinite 
length in affine space. In that setting, the existence of the 
invariant Hilbert scheme was proved by 
Haiman and Sturmfels \cite{Haiman_Sturmfels02} for 
diagonalizable groups, and (building on their work) 
in \cite{Alexeev_BrionI} for connected reductive groups. 
The connectedness assumption for $G$ is harmless and will 
be made throughout this paper.

We first present an alternative proof for the
existence of $\Hilb^G_h(V)$, by realizing it as a closed 
subscheme of $\Hilb^T_h(V//U)$; here $T\subseteq G$ is 
a maximal torus, $U\subset G$ is a maximal unipotent 
subgroup normalized by $T$, and $V//U$ is the 
categorical quotient. (Since $V//U$ is an affine 
$T$-variety, the existence of $\Hilb^T_h(V//U)$ follows
from \cite{Haiman_Sturmfels02}.) 
 
This enables us to construct another version
of the Hilbert scheme: it classifies affine $G$-schemes 
$X$ equipped with a $T$-equivariant isomorphism 
$X//U\rightarrow Y$, where $Y=\Spec(A)$ is a fixed 
affine $T$-scheme (we then say that $X$ has type $Y$). 
This moduli scheme $\M_Y$ may be described more 
concretely, in terms of the rational $G$-module $R$ 
such that the $U$-fixed subspace in $R$ 
(the span of the highest weight vectors) equals $A$:
in fact, $\M_Y$ parameterizes those algebra multiplication 
laws on $R$ that are $G$-equivariant and extend the 
multiplication of $A$.

We show that $\M_Y$ is an affine scheme of finite type,
and we equip it with an action of the adjoint torus 
$T_{\ad}$ (the quotient of $T$ by the center of $G$)
such that all orbit closures have a common point $X_0$. 
The latter corresponds to the ``most degenerate 
multiplication law'' on $R$, where the product of any two 
simple submodules is generated by the product of their 
highest weight vectors. Equivalently, $X_0$ is horospherical
in the sense of \cite{Knop90}.

These properties of $\M_Y$ follow easily from its 
description in terms of multiplication laws, except for 
the crucial fact that $\M_Y$ is of finite type. So we use
an indirect approach which relates $\M_Y$ to an invariant
Hilbert scheme: every algebra $R$ as above is generated by 
any $G$-submodule containing a set of generators for $A$, 
so that $X=\Spec(R)$ is equipped with a closed immersion into 
a fixed $G$-module $V$. And, of course, the multiplicities
of the $G$-module $R$ are fixed, since its highest weight 
vectors are.

Specifically, we define an open subscheme $\Hilb^G_h(V)_0$
of the invariant Hilbert scheme, equipped with a morphism
$f:\Hilb^G_h(V)_0 \rightarrow \Hilb^T_h(V_U)$ 
where $V_U$ is the space of co-invariants of $U$ in $V$. 
We show that $\M_Y$ identifies with the fiber of $f$ at $Y$,
for any closed $T$-stable subscheme $Y$ of $V_U$. Then
we define an action of $T_{\ad}$ on $\Hilb^G_h(V)$ that
stabilizes $\Hilb^G_h(V)_0$; it turns out that $f$ is a 
good quotient in the sense of geometric invariant theory. 
In addition, the ``most degenerate multiplication law''
yields a $T_{\ad}$-invariant section $s$ of $f$. 

The construction of $\M_Y$ gives new insight into several
methods and results concerning reductive group actions. 
For example, the degeneration of an affine $G$-variety
$X$ to a horospherical one $X_0$, which plays an 
important role in \cite{Popov86}, \cite{Vinberg86},
\cite{Knop90}, \cite{Grosshans97}, is obtained by taking
the $T_{\ad}$-orbit closure of $X$ in $\M_Y$ where 
$Y=X//U$. Consider the normalization of this orbit closure,
an affine toric variety; then a deep result of Knop 
\cite{Knop96} implies that the corresponding monoid 
(of characters) is generated by all simple roots of 
a certain root system $\Phi_X$. As a consequence, the
normalization of any $T_{\ad}$-orbit in $\M_Y$ is
isomorphic to a $T_{\ad}$-module.

Of special interest are the multiplicity-free 
$G$-varieties, i.e., those affine $G$-varieties
such that all multiplicities of their coordinate ring
are $0$ or $1$. If $G$ is a torus, then these varieties
are the (possibly non-normal) toric varieties, and 
the corresponding invariant Hilbert scheme is the
toric Hilbert scheme of Peeva and Stillman 
\cite{Peeva_Stillman02}, see also 
\cite{Haiman_Sturmfels02}. For arbitrary $G$, we show 
that every finite-dimensional $G$-module $V$ contains only 
finitely many multiplicity-free subvarieties, up to
the action of the centralizer of $G$ in $\GL(V)$. 
As a consequence, $\M_Y$ contains only finitely may
$T_{\ad}$-orbits, for any multiplicity-free $Y$.

In fact, those moduli schemes arising from reductive varieties 
(a nice class of multiplicity-free varieties, studied in 
\cite{Alexeev_BrionI}) turn out to be isomorphic to affine 
spaces, via the Vinberg family of \cite[7.5]{Alexeev_BrionI}.
This will be developed elsewhere.

\medskip

\noindent
{\sl Acknowledgements}. Part of this work was done during a staying of
the second author at the University of Georgia in January, 2003.

Both authors wish to thank the referee for his careful reading and his 
detailed and helpful comments.


\section{Existence of moduli spaces}

\label{sec: existence}


\subsection{Notation and preliminaries}

\label{subsec: notation}

We will consider algebraic groups and schemes over a fixed
algebraically closed field $k$, of characteristic zero. In addition,
all schemes will be assumed to be Noetherian.

Let $G$ be a connected reductive group, $B\subseteq G$ a Borel
subgroup with unipotent radical $U$, and $T\subseteq B$ a maximal
torus; then $B=TU$. Let $B^-=TU^-$ be the opposite Borel subgroup,
that is, $B^-\cap B = T$.

The character group $\cX(T)$ is the weight lattice
of $G$, denoted $\Lambda$ and identified with the character group
of $B$. The center $Z(G)$ is contained in $T$; the quotient $T/Z(G)$
is denoted $T_{\ad}$ and called the adjoint torus, a maximal
torus of the adjoint group $G_{\ad}=G/Z(G)$. The character group
$\cX(T_{\ad})\subseteq \Lambda$ is generated by the set 
$\Phi=\Phi(G,T)$ of roots, that is, of non-zero weights of $T$ 
in the Lie algebra $\fg$ of $G$. The choice of $B$ defines the subset 
$\Phi^+=\Phi(B,T)$ of positive roots; the corresponding subset $\Pi$ 
of simple roots is a basis of $\cX(T_{\ad})$. We write 
$\cX(T_{\ad})=\bZ\Pi$, and denote $\bN\Pi$ the submonoid generated 
by $\Pi$. This defines a partial ordering $\le$ on $\Lambda$, where 
$\mu\le\lambda$ if and only if $\lambda-\mu\in\bN\Pi$.

A rational $G$-module is a $k$-vector space $V$ equipped with 
a group homomorphism $G\rightarrow \GL(V)$ such that any $x\in V$ 
is contained in some finite-dimensional $G$-stable subspace $V_x$,
where the restriction $G\rightarrow \GL(V_x)$ is a morphism
of algebraic groups. Then $V$ is a direct sum of simple modules, and
these are finite-dimensional. 

Recall that any simple $G$-module $V$ contains a unique 
line of eigenvectors of $B$. Moreover, $V$ is uniquely determined by
the corresponding weight $\lambda\in \Lambda$; we write
$V=V(\lambda)$. Then $\lambda$ is the (unique) highest weight of $V$
(for the partial ordering $\le$), and occurs with multiplicity
$1$. The assignment $V(\lambda)\mapsto \lambda$ is a bijection from
the set of isomorphism classes of simple $G$-modules, to the subset
$\Lambda^+\subseteq \Lambda$ of dominant weights; the trivial 
one-dimensional $G$-module is $V(0)$. The dual module $V(\lambda)^*$
is simple with highest weight $-w_o\lambda$, where $w_o$ is the
longest element of the Weyl group of $(G,T)$. 

We have an isomorphism
$$
V(\lambda)^* \simeq \{f\in k[G]~\vert~ f(gb) = \lambda(b) f(g)
~\forall g\in G,~b\in B\},
$$
where $k[G]$ denotes the algebra of regular functions on $G$. 
The evaluation at the identity element $e\in G$ yields a linear 
form on $V(\lambda)^*$, eigenvector of $B$ of weight $\lambda$,
that is, a highest weight vector $v_{\lambda}\in V(\lambda)$.
This choice of highest weight vectors in simple $G$-modules is
compatible with tensor products, in the following sense: the
multiplication in $k[G]$ restricts to a surjective map of 
$G$-modules 
$V(\lambda)^*\otimes_k V(\mu)^* \rightarrow V(\lambda+\mu)^*$,
compatible with evaluations at $e$. This yields an injective map
$$
V(\lambda+\mu) \hookrightarrow V(\lambda)\otimes_k V(\mu),~
v_{\lambda+\mu}\mapsto v_{\lambda}\otimes v_{\mu}.
$$

For any rational $G$-module $M$, we have an
equivariant isomorphism
\begin{eqnarray*}
\bigoplus_{\lambda\in\Lambda^+} 
\Hom^G(V(\lambda),M)\otimes_k V(\lambda) \rightarrow M,~
\sum u_{\lambda} \otimes x_{\lambda} \mapsto 
\sum u_{\lambda}(x_{\lambda}).
\end{eqnarray*}
So the dimension of the $k$-vector space $\Hom^G(V(\lambda),M)$
is the multiplicity of $V(\lambda)$ in $M$; the image in $M$ of
$\Hom^G(V(\lambda),M)\otimes_k V(\lambda)$
is the \emph{isotypical component of $M$ of type $V(\lambda)$},
denoted $M_{(\lambda)}$. The set of all $\lambda\in\Lambda^+$ such that
$M_{(\lambda)}\ne 0$ is the \emph{weight set of} $M$, denoted
$\Lambda^+_M$. 

The multiplicity of $V(\lambda)$ in $M$ can also be read off
the subspace $M^U$ of $U$-fixed points in $M$. For $M^U$ is a
rational $T$-module, and $V(\lambda)^U=kv_{\lambda}$; thus, 
$\Hom^G(V(\lambda),M)$ is isomorphic to $M^U_{\lambda}$, 
the $\lambda$-weight space of $M^U$, via evaluation at 
the highest weight vector $v_{\lambda}$. Note that
\begin{eqnarray*}
M^U_0=M^B=M^G=M_{(0)}.
\end{eqnarray*}
We say that $M$ is \emph{multiplicity-finite}
(resp.~\emph{multiplicity-free}) if each $k$-vector space
$M_{\lambda}^U$ is finite-dimensional (resp.~of dimension at most
$1$).

An \emph{affine $G$-scheme} is an affine scheme $X=\Spec(R)$, of
finite type, equipped with an action of $G$; then $R=k[X]$ is called a
$G$-\emph{algebra}, i.e., $R$ is a finitely generated $k$-algebra 
equipped with an action of $G$ by rational automorphisms. 
We put $\Lambda^+_X=\Lambda^+_R$, the weight set of $X$.
For example, any finite-dimensional $G$-module $V$ defines an affine
$G$-scheme, with corresponding algebra $\Sym_k(V^*)$. By abuse of
notation, we still denote that subscheme by $V$.

There are obvious notions of multiplicity-finite 
(resp.~multiplicity-free) affine $G$-schemes; for simplicity, we will 
often drop the adjective ``affine''. 

As an example, since $k[G]^U_{\lambda}=V(\lambda)^*$ for all
$\lambda$, the group $G$ is multiplicity-finite (for its action on
itself by right multiplication), with weight monoid $\Lambda^+$. 
It follows that any affine $G$-variety containing a dense $G$-orbit is 
multiplicity-finite. On the other hand, the multiplicity-free
$G$-varieties are those containing a dense $B$-orbit or, 
equivalently, a dense spherical $G$-orbit.

The above definitions extend to families, in the following sense.

\begin{definition}\label{fam}
Given a scheme $S$, a \emph{family of affine $G$-schemes} over $S$ is
a scheme $\cX$ equipped with an action of $G$ and with a morphism 
$\pi:\cX \rightarrow S$, such that $\pi$ is affine, of finite type,
and $G$-invariant. 
\end{definition}

(In other words, $\cX$ is an affine $S$-scheme of finite type,
equipped with an action $a:G_S\times_S\cX \rightarrow \cX$
of the constant group scheme $G_S = G\times_{\Spec(k)} S$.)
Then the sheaf of $\cO_S$-algebras $\cR=\pi_*\cO_{\cX}$ is equipped
with a compatible $G$-action; we say that $\cR$ is a \emph{sheaf of 
$\cO_S$-$G$-algebras}. We may now formulate a basic (and well-known)
finiteness result. 

\begin{lemma}\label{dec}
With the preceding notation, we have an isomorphism of
$\cO_S$-$G$-modules:
\begin{eqnarray*}
\cR\simeq\bigoplus_{\lambda\in\Lambda^+} 
\cR_{(\lambda)},
\end{eqnarray*}
where $\cR_{(\lambda)}\simeq \cR^U_{\lambda}\otimes_k V(\lambda)$ as
$\cO_S$-$G$-modules. Moreover, both $\cR_{(0)}=\cR^G$ and $\cR^U$
are sheaves of finitely generated $\cO_S$-algebras, and every
$\cR_{\lambda}^U$ is a coherent sheaf of $\cR^G$-modules.
\end{lemma}

\begin{proof}
We may assume that $S$, and hence $\cX$, are affine. Then 
$\cR$ is a rational $G$-module by \cite[I.1]{Mumford_GIT3ed}, and 
the first assertion follows from the structure of these
modules.

Let $R=\Gamma(\cX,\cO_{\cX})$ and $A=\Gamma(S,\cO_S)$. Then $R$ is a
finitely generated $A$-$G$-algebra. Thus, it is generated by a
finite-dimensional $G$-submodule $V\subseteq R$. Then the
$k$-algebra $\Sym_k(V)^G$ is finitely generated, by the
Hilbert-Nagata theorem. Since the $A$-algebra $R^G$ is a quotient 
of $A\otimes_k \Sym_k(V)^G$, it is finitely generated. 
Likewise, each $R_{\lambda}^U$ is a finitely generated
$A$-module, and $R^U$ is a finitely generated $A$-algebra, see
e.g. \cite[Theorems 5.6, 9.1]{Grosshans97}.
\end{proof}
Next put $\cX//G=\Spec_S(\cR^G)$ and let 
$$
p=p_{\cX,G}:\cX\rightarrow \cX//G
$$ 
be the morphism corresponding to the inclusion 
$\cR^G \subseteq \cR$. Then $p$ is a \emph{good quotient} by 
$G$, that is: $p$ is affine and $G$-invariant, and induces an
isomorphism $\cO_{\cX//G}\rightarrow (p_*\cO_{\cX})^G$. It follows
e.g. that $p$ is surjective, and sends any closed $G$-stable subset
to a closed subset; as a consequence, $p$ is universally
submersive (for these facts, see \cite[I.2]{Mumford_GIT3ed}). 
Moreover, $\pi = (\pi//G) \circ p$, where  
$$
\pi//G:\cX//G\rightarrow S
$$ 
is an affine morphism of finite type.

Likewise, we obtain a factorization $\pi = (\pi//U) \circ p_{\cX,U}$,
where 
$$
p_{\cX,U}:\cX\rightarrow \cX//U
$$ 
is a good quotient by $U$, and 
$$
\pi//U:\cX//U\rightarrow S
$$ 
is a family of affine $T$-schemes. Moreover, $(\cX//U)//T =\cX//B$ is
mapped isomorphically to $\cX//G$. 

Many properties of an affine $G$-scheme $X$ can be read off $X//U$,
e.g., $X$ is reduced (resp. a variety, normal) if and only if
$X//U$ is, see \cite{Popov86}. Since $\Lambda^+_X = \Lambda^+_{X//U}$,
it follows that $\Lambda^+_X$ is a finitely generated monoid, 
for any affine $G$-variety $X$.

Given two families of affine $G$-schemes $\pi:\cX\rightarrow S$
and $\pi':\cX'\rightarrow S$, we denote by $\Hom^G_S(\cX,\cX')$
the set of $G$-equivariant morphisms $\cX\rightarrow\cX'$ over $S$;
these will just be called morphisms. Any morphism
$\gamma:\cX\rightarrow \cX'$ yields a $T$-equivariant morphism
$\gamma//U:\cX//U\rightarrow \cX'//U$ over $S$. 

With these notations, we may record another easy finiteness result.

\begin{lemma}\label{aut}
(i) The map $\Hom^G_S(\cX,\cX')\rightarrow\Hom^T_S(\cX//U,\cX'//U)$
is injective.

\noindent
(ii) For any two multiplicity-finite $G$-schemes $X$, $X'$, 
the contravariant functor (Schemes) $\rightarrow$ (Sets),
$S\mapsto \Hom^G_S(X\times S,X'\times S)$ is representable by 
an affine scheme of finite type denoted $\Hom^G(X,X')$.

Moreover, the group functor $S\mapsto \Aut^G_S(X\times S)$ is 
representable by a linear algebraic group $\Aut^G(X)$, open in 
$\Hom^G(X,X)$. In particular, the Lie algebra of $\Aut^G(X)$ is 
the space $\Der^G(X)$ of $G$-equivariant derivations of $X$.

If, in addition, $X$ is multiplicity-free, then $\Aut^G(X)$ 
is diagonalizable.
\end{lemma}

\begin{proof}
(i) follows at once from the corresponding statement for 
rational $G$-modules. 

(ii) Choose a $G$-equivariant closed immersion 
$X'\hookrightarrow V$, where $V$ is a finite-dimensional 
$G$-module. Let $R=\Gamma(X,\cO_X)$ and $A=\Gamma(S,\cO_S)$. 
Then
$$
\Hom^G_S(X\times S,V\times S) = 
\Hom^G_A(A\otimes_k\Sym_k(V^*),A\otimes_k R) \simeq
A\otimes_k (R\otimes_k V)^G.
$$
So $S\mapsto \Hom^G_S(X\times S,V\times S)$ is represented 
by the affine space $(R\otimes_k V)^G$. The latter is 
finite-dimensional, since $X$ is multiplicity-finite. 
Moreover, denoting by $I$ the ideal of $X'$ in $\Sym_k(V^*)$, 
the restriction to $I$ yields a morphism
$$
\varphi: (R\otimes_k V)^G = \Hom^{G-\alg}(\Sym_k(V^*),R) 
\rightarrow \Hom^{G-\mod}(I,R),
$$
with evident notation. Further, $\Hom^{G-\mod}(I,R)$ is an affine
space, possibly infinite-dimensional.
Clearly, the scheme-theoretic fiber $\varphi^{-1}(0)$ 
represents $S\mapsto \Hom^G_S(X\times S,X'\times S)$.

In particular, $\Hom^G(X,X)$ is a closed subscheme of the affine space
$$
\End^G(\bigoplus_{\lambda\in F} R_{(\lambda)})=
\End^G(\bigoplus_{\lambda\in F} R_{\lambda}^U\otimes_k V(\lambda))=
\prod_{\lambda\in F} \End(R_{\lambda}^U), 
$$
where $F$ is any finite subset of $\Lambda^+$ such that the algebra
$R$ is generated by its isotypical components $R_{(\lambda)}$,
$\lambda\in F$. By Lemma \ref{dec}, $\Aut^G(X)$ is the intersection of 
$\Hom^G(X,X)$ with the open subset 
$\prod_{\lambda\in F} \GL(R_{\lambda}^U)\subset
\prod_{\lambda\in F} \End(R_{\lambda}^U)$. So $\Aut^G(X)$ is a linear
algebraic group; arguing as above, we see that it represents
$S\mapsto \Aut^G_S(X\times S)$.

Finally, if $X$ is multiplicity-free, then each non-zero
$R_{\lambda}^U$ is a line, so that $\Aut^G(X)$ is diagonalizable.
\end{proof}


\subsection{The invariant Hilbert scheme}

\label{subsec: invariant}

Consider a family $\pi:\cX\rightarrow S$ of affine $G$-schemes and
put $\cR=\pi_*\cO_{\cX}$. By Lemma \ref{dec}, the morphism $\pi$ is 
flat if and only if $\pi//U$ is flat; then $\pi//G$ is flat as well. 
If, in addition, $\pi//G$ is finite, then every $\cR_{\lambda}^U$ 
is a flat, coherent sheaf of $\cO_S$-modules, and hence locally
free of finite rank. This motivates the following

\begin{definition}\label{flt}
Given a function $h:\Lambda^+\to \bN$, a family of 
multiplicity-finite $G$-schemes $\pi:\cX\rightarrow S$ 
\emph{has Hilbert function $h$}, if every sheaf of $\cO_S$-modules 
$\cR_{\lambda}^U$ is locally free of constant rank $h(\lambda)$.
\end{definition}

Then $\pi$ is flat, and $\pi//U:\cX//U\rightarrow S$ is a family of
multiplicity-finite $T$-schemes with the same Hilbert function
$h$. Note that $h(0)\ge 1$, so that $\pi//G$ is finite and
surjective. As a consequence, $\pi=(\pi//G)\circ p$ is surjective and
maps closed $G$-stable subsets to closed subsets.

\begin{definition}\label{hil}
Given an affine $G$-scheme $X$ and a function $h:\Lambda^+\to\bN$,
the \emph{Hilbert functor} is the contravariant functor
$\cHilb^G_h(X)$: (Schemes) $\rightarrow$ (Sets) assigning to any
scheme $S$ the set of closed $G$-stable subschemes 
$\cX\subseteq X\times S$ such that the projection $\pi:\cX\to S$ is a
family of multiplicity-finite $G$-schemes with Hilbert function
$h$.
\end{definition}

In the case where $G=T$ is a torus and $X$ is a finite-dimensional
$T$-module,  the functor $\cHilb^T_h(X)$ is represented by a
quasi-projective scheme $\Hilb^T_h(X)$, as follows from 
\cite[Theorem 1.1]{Haiman_Sturmfels02}. We will need the following 
straightforward consequence of this result.

\begin{lemma}\label{rep0}
For any affine $T$-scheme $Y$ and any function $h:\Lambda^+\to\bN$,
the functor $\cHilb^T_h(Y)$ is represented by a quasi-projective
scheme $\Hilb^T_h(Y)$.
\end{lemma}

\begin{proof}
Choose a $T$-equivariant closed immersion $Y\hookrightarrow E$, where
$E$ is a finite-dimensional $T$-module; let $I\subset \Sym_k(E^*)$
be the corresponding ideal. Then any family
$\pi:\cX\rightarrow S$ in $\cHilb^T_h(Y)(S)$ yields a family in
$\cHilb^T_h(E)(S)$, that is, a morphism $S\rightarrow \Hilb^T_h(E)$
such that $\cX$ is the pull-back of the universal family 
$$
\Univ^T_h(E)\subseteq E\times \cHilb^T_h(E).
$$ 
Let $\pi:\Univ^T_h(E)\rightarrow\Hilb^T_h(E)$ be the projection.
Then each eigenspace 
$$
\cF_{\lambda}=(\pi_*\cO_{\Univ^T_h(E)})_{\lambda}
$$ 
is a locally free sheaf of rank $h(\lambda)$ on $\Hilb^G_h(E)$,
equipped with a map 
$$
q_{\lambda}:\Sym_k(E^*)_{\lambda} \rightarrow
\Gamma(\Hilb^G_h(E),\cF_{\lambda}). 
$$
Moreover, the image of $q_{\lambda}(I_{\lambda})$ under the pull-back 
$$
\Gamma(\Hilb^G_h(E),\cF_{\lambda}) \rightarrow
\Gamma(S,\cR_{\lambda})
$$ 
is zero, since $\cX\subseteq Y\times S$. The intersection over all 
$\lambda$ of the zero subschemes of the spaces 
$q_{\lambda}(I_{\lambda})$ is a closed subscheme of $\Hilb^G_h(E)$, 
that clearly represents $\Hilb^G_h(Y)$.
\end{proof}

We now come to our first main result. 

\begin{theorem}\label{rep1}
For any affine $G$-scheme $X$ and any function $h:\Lambda^+\to\bN$,
the functor $\cHilb^G_h(X)$ is represented by a closed subscheme 
$\Hilb^G_h(X)$ of $\Hilb^T_h(X//U)$.
\end{theorem}

\begin{proof}
Let $S$ be a scheme. Then $\cHilb^G_h(X)(S)$ is the set of
those $G$-stable ideal sheaves 
$\cI\subseteq \cO_S \otimes_k \Gamma(X,\cO_X)$
such that each sheaf
$(\cO_S \otimes_k \Gamma(X,\cO_X)/\cI)_{\lambda}^U$ is locally free of
rank $h(\lambda)$. Such an ideal sheaf $\cI$ is uniquely determined by
$\cI^U$; the latter is a $T$-stable ideal sheaf of 
$\cO_S\otimes_k \Gamma(X//U,\cO_{X//U})$, which yields a point of 
$\cHilb^T_h(X//U)(S)=\Hom(S,\Hilb^T_h(X//U))$.

Conversely, a given $\cJ\in \Hom(S,\Hilb^T_h(X//U))$ equals $\cI^U$
for some $\cI\in\cHilb^G_h(X)(S)$, if and only if the span
\begin{eqnarray*}
\langle G\cdot \cJ\rangle\subseteq \cO_S \otimes_k \Gamma(X,\cO_X)
\end{eqnarray*} 
of all $G$-translates of $\cJ$, is a sheaf of ideals. Since 
$\langle G\cdot \cJ\rangle$ is a sheaf of $\cO_S$-$G$-modules, this 
amounts to 
\begin{eqnarray*}
(\langle G\cdot \cJ\rangle\cdot \Gamma(X,\cO_X))^U\subseteq \cJ.
\end{eqnarray*}
To express this condition in more concrete terms, consider three
dominant weights $\lambda$, $\mu$, $\nu$, and a $B$-eigenvector
$v\in V(\lambda)\otimes_k V(\mu)$, of weight $\nu$. We may write
$$
v=\sum_{i\in F} a_i(g_i\cdot v_{\lambda})\otimes (h_i \cdot v_{\mu}),
$$ 
where $F$ is a finite set, $a_i\in k$, and $g_i, h_i\in G$.
Then, for any local section $\varphi$ of $\cJ_{\lambda}$ and for any 
$\psi\in\Gamma(X,\cO_X)^U_{\mu}$, we obtain a local section 
$\sum_{i\in F} a_i (g_i\cdot \varphi) (h_i\cdot \psi)$
of $\cO_S\otimes_k \Gamma(X,\cO_X)^U_{\nu}$, and hence a morphism
$$
\cJ_{\lambda}\otimes_k \Gamma(X,\cO_X)^U_{\mu} \rightarrow 
\cO_S\otimes_k \Gamma(X,\cO_X)^U_{\nu}
$$ 
of sheaves over $\cO_S$. Composing with the quotient by $\cJ_{\nu}$
yields a morphism of $\cO_S$-sheaves
\begin{eqnarray*}
\cJ_{\lambda}\otimes_k \Gamma(X//U,\cO_{X//U})_{\mu} \rightarrow
\cO_S\otimes_k \Gamma(X//U,\cO_{X//U})_{\nu}/\cJ_{\nu}.
\end{eqnarray*}
Now our condition is the vanishing of all such morphisms. The latter
may be regarded as sections of certain locally free sheaves over
$\Hilb^T_h(X//U)$, arising from its universal family. So their
vanishing defines a closed subscheme of $\Hilb^T_h(X//U)$, which
represents $\cHilb^G_h(X)$.
\end{proof}

We say that $\Hilb^G_h(V)$ is the invariant Hilbert scheme.
Assigning to any family its quotient by $G$, we obtain
a morphism 
$$
\eta:\Hilb^G_h(V)\rightarrow \Hilb_{h(0)}(V//G)
$$
to the punctual Hilbert scheme that parameterizes closed subschemes
of length $h(0)$ in $V//G$. One may check that $\eta$ is proper;
it generalizes the Nakamura morphism for finite groups, see e.g.
\cite{Nakamura01}.


\subsection{Moduli $\M_Y$ of affine schemes of type $Y$}

\label{subsec: moduli}
To construct this moduli scheme, we need a relative version 
of the invariant Hilbert scheme. Consider two affine $G$-schemes
$X$, $Y$, an equivariant morphism $f:X\rightarrow Y$, and a function
$h:\Lambda^+ \rightarrow \bN$. Let $\cHilb^G_h(f)$ be the functor
assigning to each scheme $S$ the subset of $\cHilb^G_h(X)(S)$
consisting of those $\cX\subseteq X\times S$ such that the morphism 
$$
f\times \id: \cX\rightarrow Y\times S
$$ 
is a closed immersion.

\begin{lemma}\label{rel}
The functor $\cHilb^G_h(f)$ is represented by an open subscheme
$\Hilb^G_h(f)$ of $\Hilb^G_h(X)$, equipped with a morphism 
$\Hilb^G_h(f) \rightarrow \Hilb^G_h(Y)$.
\end{lemma}

\begin{proof}
Consider a family $\pi:\cX\rightarrow S$ in $\cHilb^G_h(X)(S)$; put 
$\cR=\pi_*\cO_{\cX}$. By Theorem \ref{rep1}, we have a morphism 
$S \rightarrow \Hilb^T_h(X//U)$ such that the family
$\pi//U: \cX//U \rightarrow S$ is the pull-back of the universal
family $\Univ^T_h(X//U) \rightarrow \Hilb^T_h(X//U)$. As a
consequence, there exists a finite subset $F\subset\Lambda^+$ 
(depending only on $X$ and $h$)
such that the sheaf of $\cO_S$-algebras $\cR^U$ is
generated by its weight spaces $\cR^U_{\lambda}$, $\lambda\in F$. 

Now $\pi$ lies in $\cHilb^G_h(f)(S)$ if and only if the map 
$$
(f\times\id)^{\#}:\cO_S\otimes_k\Gamma(Y,\cO_Y)\rightarrow \cR
$$ 
is surjective; equivalently, its restriction
\begin{eqnarray*}
\cO_S\otimes_k\Gamma(Y,\cO_Y)^U_{\lambda} \rightarrow \cR^U_{\lambda}
\end{eqnarray*}
is surjective for all $\lambda\in F$. This yields finitely many open
conditions on $\Hilb^G_h(X)\subseteq \Hilb^T_h(X//U)$, and hence an
open subset $\Hilb^G_h(f)$ that represents
$\cHilb^G_h(f)$. Moreover, assigning to $\cX$ its image in $Y\times S$
yields a natural transformation 
$\cHilb^G_h(f)\rightarrow \cHilb^G_h(Y)$, and hence the desired 
morphism.
\end{proof}

Next let $V$ be a finite-dimensional $G$-module and let $V_U$ be the
space of \emph{co-invariants} of $U$, i.e., $V_U$ is the quotient of $V$ 
by the span of the elements $u\cdot v - v$, where $u\in U$ and
$v\in V$. In other words, the quotient map 
$$
q:V\rightarrow V_U
$$ 
is the universal $U$-invariant linear map with source $V$. Since $T$
normalizes $U$, it acts on $V_U$, and $q$ is equivariant. 

For any $\lambda\in\Lambda^+$, the $T$-module $V(\lambda)_U$ is
one-dimensional with weight $w_o\lambda$, the lowest weight
of $V(\lambda)$. Moreover, $q$ maps isomorphically $V(\lambda)^{U^-}$
to $V(\lambda)_U$. It follows that the restriction $V^{U^-}\to V_U$ is
an isomorphism. Thus, the category of finite-dimensional $G$-modules
is equivalent to that of finite-dimensional $T$-modules with weights
in $-\Lambda^+$, via $V\mapsto V_U$.

The map $q:V\rightarrow V_U$ factors as $f\circ p$, where
$p=p_{V,U}:V\rightarrow V//U$  is the quotient, and 
\begin{eqnarray*}
f:V//U\rightarrow V_U
\end{eqnarray*}
is a surjective, $T$-equivariant morphism.

Let $\cHilb^G_h(V)_0$ be the functor assigning to every scheme $S$
the subset of $\cHilb^G_h(V)(S)$ consisting of those subfamilies
$\cX\subseteq V\times S$ such that the morphism 
$$
f//U\times \id: \cX//U\rightarrow V_U\times S
$$
is a closed immersion. Then Theorem \ref{rep1} and Lemma \ref{rel}
immediately imply the following result.

\begin{theorem}\label{rep2}
The functor $\cHilb^G_h(V)_0$ is represented by an open subscheme
$\Hilb^G_h(V)_0$ of $\Hilb^G_h(V)$, equipped with a morphism
\begin{eqnarray*}
f:\Hilb^G_h(V)_0 \rightarrow \Hilb^T_h(V_U).
\end{eqnarray*}
\end{theorem}

\begin{definition}\label{typ}
Given an affine $T$-scheme $Y$, a 
\emph{family of affine $G$-schemes of type $Y$} over a scheme $S$
consists of a family of affine $G$-schemes $\pi:\cX\rightarrow S$,
together with an isomorphism 
$\varphi:\cX//U\rightarrow Y\times S$ of families of affine 
$T$-schemes over $S$.

We say that the families of type $Y$
$(\pi:\cX\rightarrow S,~\varphi:\cX//U\rightarrow Y\times S)$ 
and $(\pi':\cX'\rightarrow S,~\varphi':\cX'//U\rightarrow Y\times S)$
are \emph{equivalent} if there exists a morphism
$\psi:\cX\rightarrow\cX'$ such that
$\varphi'\circ\psi//U=\varphi$.
\end{definition}

With the preceding notation, $\pi$ is flat, since the family $\pi//U$ 
is trivial. Moreover, $\pi$ is a family of multiplicity-finite 
$G$-schemes, if and only if $Y$ is multiplicity-finite; then $\pi$ 
and $Y$ have the same Hilbert function. On the other hand, $\psi$ is 
an isomorphism, since $\psi//U$ is. Likewise, any self-equivalence 
is the identity.

\begin{definition}\label{mod}
Let $Y$ be a multiplicity-finite $T$-scheme. The 
\emph{moduli functor of affine $G$-schemes of type $Y$}
is the contravariant functor 
$\cM_Y$: (Schemes) $\rightarrow$ (Sets) 
assigning to any $S$ the set of equivalence classes of families of
affine $G$-schemes of type $Y$ over $S$. 
\end{definition}

By Lemma \ref{aut}, the group $\Aut^T_S(Y\times S)$ equals
$\Hom(S,\Aut^T(Y))$. It acts naturally on $\cM_Y(S)$, and the
isotropy group of the equivalence class of 
$(\pi:\cX\rightarrow S,~\varphi:\cX//U\rightarrow Y\times S)$
is the image of the map
$$
\Aut^G_S(\cX) \rightarrow \Aut^T_S(\cX//U) \simeq
\Aut^T_S(Y\times S) = \Hom(S,\Aut^T(Y))
$$ 
(which is injective, by Lemma \ref{aut}). The orbit space
$\cM_Y(S)/\Aut^T_S(Y\times S)$ is the set of isomorphism classes
of families $\cX$ over $S$, such that $\cX//U$ is isomorphic to
the trivial family $Y\times S$.

We may now describe $\cM_Y$ in terms of invariant Hilbert schemes.

\begin{theorem}\label{rep3}
Let $Y$ be a multiplicity-finite $T$-scheme, with Hilbert
function $h$. Choose a finite-dimensional $T$-module $E$ such that $Y$ 
admits a closed $T$-equivariant immersion into $E$. Let $V$ be the 
$G$-module such that $V_U=E$, and let $\M_Y$ be the fiber at $Y$ of
the morphism $f:\Hilb^G_h(V)_0\rightarrow \Hilb^T_h(V_U)$. Then 
$\M_Y$ represents the functor $\cM_Y$. 

In particular, $\M_Y$ is independent of $E$ and admits an action of
$\Aut^T(Y)$; the orbits are in bijection with the isomorphism classes
of affine $G$-schemes $X$ such that $X//U\simeq Y$.
\end{theorem}

\begin{proof}
Let $(\pi:\cX\rightarrow S,~\varphi:\cX//U\rightarrow Y\times S)$ be a
family of affine $G$-schemes of type $Y$; let
$\cR=\pi_*\cO_{\cX}$. Then we have a morphism of $T$-modules 
$E^*\rightarrow \Gamma(S,\cR^U)$ that extends uniquely to a morphism
of $G$-modules $V^*\rightarrow \Gamma(S,\cR)$. This yields a morphism
of sheaves of $\cO_S$-$G$-algebras 
$$
\cO_S\otimes_k\Sym_k(V^*)\rightarrow \cR
$$ 
lifting $\cO_S\otimes_k\Sym_k(E^*)\rightarrow \cR^U$. Since the latter
is surjective, we obtain a closed $G$-equivariant immersion
$\cX\hookrightarrow V\times S$ over $S$, such that the composition 
\begin{eqnarray*}
\cX//U\rightarrow V//U\times S\rightarrow V_U\times S = E\times S
\end{eqnarray*}
is nothing but $\varphi:\cX//U\rightarrow Y\times S$ followed by 
the inclusion $Y\times S\subseteq E\times S$. So, by Theorem
\ref{rep2}, the image of $\cX$ in $V\times S$ yields an $S$-point of
the fiber $f^{-1}(Y)$. Clearly, this point only depends on the
equivalence class of the family $(\pi,\varphi)$.

For the converse, note that an $S$-point of $\Hilb^G_h(V)$ is nothing
but an equivalence class of pairs $(\pi,\iota)$, where
$\pi:\cX\rightarrow S$ is a family of multiplicity-finite
$G$-schemes with Hilbert function $h$, and 
$\iota:\cX\rightarrow V\times S$ is a closed $G$-equivariant
immersion over $S$. Thus, any $S$-point of $f^{-1}(Y)$ yields an
equivalence class of pairs as above, with the additional assumption
that the composition 
$\cX//U\rightarrow V//U\times S\rightarrow V_U\times S = E\times S$ 
is a closed immersion with image $Y\times S$. This yields in turn an
equivalence class of families of affine $G$-schemes of type $Y$ over
$S$.
\end{proof}

Now consider the following moduli stack (see 
\cite{Laumon_Moret-Bailly00} for definitions and basic properties 
of stacks). For every scheme 
$S$, let $M'_Y(S)$ be the category whose objects are families of 
affine $G$-schemes $\pi:X\to S$ such that for some \'etale cover 
$\{S_i\to S\}$ the schemes $(X//U)\times_S S_i$ are 
isomorphic to $Y\times_S S_i$ as $T$-schemes over 
$S_i$. Morphisms in $M'_Y(S)$ are just isomorphisms of families
of affine $G$-schemes. It is immediate that $M'_Y$ is
the quotient stack of $M_Y$ by the linear algebraic group 
$\Aut^T(Y)$. It is certainly not a Deligne-Mumford stack in 
general, since the stabilizers may be infinite. However, it is an 
algebraic Artin stack, and it is smooth if the scheme $M_Y$ is.


\subsection{The tangent space of the invariant Hilbert scheme}

\label{subsec: tangent}

We begin with a description of the tangent space of $\Hilb^G_h(V)$, 
completely analogous to that of the classical Hilbert scheme; the 
tangent space of $\M_Y$ will be determined in 
\ref{subsec: structure} below.

\begin{proposition}\label{tgt1}
Consider a finite-dimensional $G$-module $V$, a function
$h:\Lambda^+\to\bN$, and a closed point $X\in\Hilb^G_h(V)$, that is,
a $G$-stable subscheme $X\subseteq V$ with Hilbert function $h$. Let
$I\subseteq \Sym_k(V^*)$ be the ideal of $X$, and let
$R=\Sym_k(V^*)/I$. Then $\Hom_R(I/I^2,R)$ is a multiplicity-finite
$G$-module, and the Zariski tangent space $T_X \Hilb^G_h(V)$
is canonically isomorphic to $\Hom^G_R(I/I^2,R)$. 

Moreover, the space $T^1(X)$ of infinitesimal deformations of $X$ is
also a multiplicity-finite $G$-module, and we have an exact sequence
of finite-dimensional $k$-vector spaces
\begin{eqnarray*}
0 \rightarrow \Der^G(X) \rightarrow
\Hom^G(X,V) \rightarrow T_X \Hilb^G_h(V) \rightarrow T^1(X)^G
\rightarrow 0.  
\end{eqnarray*}

\end{proposition}

\begin{proof}
Let $D=k[t]/(t^2)=k\oplus k\varepsilon$, where $\varepsilon^2=0$. Then
$T_X \Hilb^G_h(V)$ is the fiber at $X$ of the map
$\Hilb^G_h(V)(\Spec D) \rightarrow \Hilb^G_h(V)(\Spec k)$. In other
words, $T_X \Hilb^G_h(V)$ consists of those $G$-stable ideals
$J\subseteq D\otimes_k\Sym_k(V^*)$ such that: 
$D\otimes_k\Sym_k(V^*)/(J,\varepsilon)=R$, and
$D\otimes_k\Sym_k(V^*)/J$ is flat over $D$. Equivalently, 
$J\cap\Sym_k(V^*)=I$ and 
$J\cap\varepsilon\Sym_k(V^*)=\varepsilon I$. For any 
$u+\varepsilon v\in J$, it follows that: $u\in I$, and the class
$v+I\in R$ is uniquely determined by $u$. This yields a map
$\varphi:I\rightarrow R,~u\mapsto v+I$. One easily checks that 
$\varphi$ lies in $\Hom^G_{\Sym_k(V^*)}(I,R)=\Hom^G_R(I/I^2,R)$, and
that $J$ is the preimage in $I\oplus \varepsilon \Sym_k(V^*)$ of the
graph of $\varphi$ in $I\oplus\varepsilon R$. Moreover, any 
$\varphi\in\Hom^G_R(I/I^2,R)$ arises from an ideal 
$J\in  T_X \Hilb^G_h(V)$. This proves the first assertion.

Next recall the exact sequence of K\"ahler differentials (over $k$)
$$
I/I^2 \rightarrow \Omega^1_{\Sym_k(V^*)}\otimes_{\Sym_k(V^*)}R
\rightarrow \Omega^1_R\rightarrow 0,
$$
that is,
$$
I/I^2 \rightarrow R\otimes_k V^* \rightarrow \Omega^1_R\rightarrow 0.
$$
Applying $\Hom_R(-,R)$, we obtain a exact sequence of
$R$-$G$-modules
$$
0 \rightarrow \Der(R) \rightarrow R\otimes_k V \rightarrow
\Hom_R(I/I^2,R) \rightarrow T^1(X) \rightarrow 0,
$$
see e.g. \cite[Exercise II.9.8]{Hartshorne77}. Since 
$R\otimes_k V = \Hom(X,V)$, this yields our exact sequence by taking 
$G$-invariants.

Since $I/I^2$ is a finitely generated $R$-$G$-module, it is
generated as an $R$-module by a finite-dimensional $G$-submodule
$E$. Then $\Hom_R(I/I^2,R)$ identifies to a $R$-$G$-submodule of
$R\otimes_k E^*$, and the latter is multiplicity-finite as a
$G$-module, since $R$ is. So $\Hom_R(I/I^2,R)$, and hence 
$T^1(X)$, are multiplicity-finite $G$-modules; in particular,
$\Hom_R^G(I/I^2,R)$ and $T^1(X)^G$ are finite-dimensional.  
The same holds for $\Hom^G(X,V)=(R\otimes_k V)^G$. Thus, 
$\Der^G(X)$ is finite-dimensional as well (this also follows 
from Lemma \ref{aut}).
\end{proof}

Denote by $\Aut^G(V)$ the automorphism group of the $G$-variety
$V$. This group acts on $\Hilb^G_h(V)$, 
and we may think of $T^1(X)^G$ as the normal space at $X$ 
to the orbit $\Aut^G(V)\cdot X$. Indeed, we may regard 
$\Der^G(V)$ as the Lie algebra of the (possibly 
infinite-dimensional) algebraic group $\Aut^G(V)$. Further, 
$\Der^G(V) = (\Sym_k(V^*) \otimes_k V)^G$
maps surjectively to $(R\otimes_k V)^G = \Hom^G(X,V)$,
so that the image of $\Hom^G(X,V)$ in $T_X \Hilb^G_h(V)$ is
$\Der^G(V)\cdot X$, the tangent space to the orbit.
Under additional assumptions, we will see that $\Aut^G(V)$
may be replaced with the centralizer $\GL(V)^G$. The latter 
is isomorphic to the product of the $\GL(V^U_{\lambda})$'s, 
and hence is a connected reductive algebraic group.

Let $X\subseteq V$ be a multiplicity-free subvariety,
that is, the closure of a spherical orbit.
Let $\cS=\Lambda^+_X$ be the corresponding weight monoid, 
then $h(\lambda)=1$ if $\lambda\in\cS$, and $h(\lambda)=0$ otherwise.
Thus, the data of $\cS$ and of $h$ are equivalent; we will write
$\Hilb^G_{\cS}(V)$ for the invariant Hilbert scheme $\Hilb^G_h(V)$.
The group $\GL(V)^G$ acts on that scheme; thus, we can consider 
the orbit $\GL(V)^G\cdot X$, where $X$ is regarded as a closed point
of $\Hilb^G_{\cS}(V)$. 

\begin{definition}\label{ndg}
Given a finite-dimensional $G$-module $V$ and a closed $G$-stable 
subvariety $X$, we say that $X$ is \emph{non-degenerate}, 
if its projections to the isotypical components of $V$ are all 
non-zero. 
\end{definition}

(In particular, the non-degenerate subvarieties of a
multiplicity-free $G$-module $V$ are those that span $V$.) 
Now Proposition \ref{tgt1} may be refined as follows.

\begin{proposition}\label{tgt2}
Let $X$ be the closure of a spherical orbit $G\cdot x$ in a 
finite-dimensional $G$-module $V$. 

\noindent
(i) If $X$ is non-degenerate, then the normal space at $X$ to 
$\GL(V)^G\cdot X$ is isomorphic to $T^1(X)^G$.

\noindent
(ii) If, in addition, the isotropy group $G_x$ is reductive, 
then $T^1(X)^G=0$. In other words, $\GL(V)^G\cdot X$ is open 
in $\Hilb^G_{\cS}(V)$.

\noindent
(iii) On the other hand, if $X$ is normal and if the boundary
$X - G\cdot x$ has codimension at least $2$ in $X$, then 
the exact sequence of Proposition \ref{tgt1} identifies with
$$
0 \rightarrow (\fg/\fg_x)^{G_x} \rightarrow V^{G_x}
\rightarrow  (V/\fg\cdot x)^{G_x} \rightarrow T^1(X)^G 
\rightarrow 0,
$$
where $\fg_x$ denotes the isotropy Lie algebra of $x$. 
As a consequence, we have an exact sequence
$$
0 \rightarrow T^1(X)^G \rightarrow H^1(G_x,\fg/\fg_x)
\rightarrow H^1(G_x,V).
$$
\end{proposition}

\begin{proof}
(i) The Lie algebra of $\GL(V)^G$ is $\End^G(V)$, so that 
the tangent space to $\GL(V)^G\cdot X$ at $X$ equals 
$\End^G(V)\cdot X$. With the notation of Proposition
\ref{tgt1}, this subspace of $\Hom^G_R(I/I^2,R)$ 
is the image of the composition
$$\CD
\End^G(V) @>{r_X}>> \Hom^G(X,V) @>>> \Hom^G_R(I/I^2,R),
\endCD
$$
where $r_X$ denotes the restriction map. By the exact
sequence of that proposition, it suffices to show 
that $r_X$ is surjective. And since the restriction
$\Hom^G(X,V) \rightarrow \Hom^G(G\cdot x,V) = V^{G_x}$
is injective, it suffices in turn to check the surjectivity
of the composition 
$$
\End^G(V)\rightarrow V^{G_x},~f \mapsto f(x).
$$

Write 
$$
V=\bigoplus_{\lambda\in F} M_{\lambda}\otimes_k V(\lambda),
$$
where $F$ is a finite subset of $\Lambda^+$, and 
$M_{\lambda}= \Hom^G(V(\lambda),V)$. Then 
$$
\End^G(V) \simeq \prod_{\lambda\in F} \End(M_{\lambda}).
$$
On the other hand, we have $\dim V(\lambda)^{G_x}\le 1$ 
for any $\lambda\in\Lambda^+$, since $G_x$ is a spherical 
subgroup of $G$. And since $X$ has a nonzero projection 
on each isotypical component
$V_{(\lambda)} = M_{\lambda}\otimes V(\lambda)$, we have
$$
x = \sum_{\lambda\in F} m_{\lambda} \otimes x_{\lambda},
$$
where each $m_{\lambda}$ is nonzero, and each $x_{\lambda}$ spans 
$V(\lambda)^{G_x}$. It follows that
$$
\End^G(V)(x) = 
\sum_{\lambda\in F} M_{\lambda}\otimes V(\lambda)^{G_x} 
= V^{G_x},
$$
as desired.

(ii) Since $G\cdot x$ is dense in $X$, the restriction
map
$\Hom_R(I/I^2,R) \rightarrow H^0(G\cdot x,\cN_{G\cdot x})$
is injective, where $\cN_{G\cdot x}$ denotes the normal 
sheaf. The latter is the $G$-linearized sheaf on $G/G_x$ 
associated with the $G_x$-module $V/\fg\cdot x$. 
This yields an injection
$$
T_X\Hilb^G_{\cS}(V) = \Hom_R^G(I/I^2,R) \rightarrow 
H^0(G\cdot x,\cN_{G\cdot x})^G = (V/\fg\cdot x)^{G_x}.
$$
But since $G_x$ is reductive, the natural map
$V^{G_x} \rightarrow (V/\fg\cdot x)^{G_x}$
is surjective. Together with the proof of (i), it 
follows that $\End^G(V)X$ equals $T_X\Hilb^G_{\cS}(V)$.

(iii) Since $X=\Spec(R)$ is normal, the $R$-modules
$$
\Der(X)=\Hom_R(\Omega^1_R,R),~\Hom(X,V)=R\otimes_k V,
~\Hom_R(I/I^2,R)
$$ 
are all reflexive. And since $\codim_X(X - G\cdot x)\ge 2$, 
the exact sequence
$$
0 \rightarrow \Der(X) \rightarrow \Hom(X,V) \rightarrow
\Hom_R(I/I^2,R)
$$
identifies with
$$
0 \rightarrow H^0(G\cdot x,\cT_{G\cdot x}) \rightarrow
H^0(G\cdot x, \cO_{G\cdot x}\otimes_k V) \rightarrow
H^0(G\cdot x,\cN_{G\cdot x}),
$$
where $\cT_{G\cdot x}$ denotes the tangent sheaf. Further, 
$\cT_{G\cdot x}$ (resp. $\cO_{G\cdot x}\otimes_k V$) 
is the $G$-linearized sheaf on $G/G_x$ associated with 
the $G_x$-module $\fg/\fg_x$ (resp. $V$). Now taking 
$G$-invariants in the former exact sequence yields our 
assertion.
\end{proof}

In the case where $G=T$ is a torus, the multiplicity-free 
varieties are those containing a dense orbit; these are 
the (possibly non-normal) affine toric varieties. They are 
classified by finitely generated submonoids of 
$\Lambda^+=\Lambda$, via $X\mapsto \Lambda_X$ and 
$\cS\mapsto \Spec k[\cS]$. Note that the map 
$T\mapsto \Aut^T(X)$ is surjective with kernel 
the intersection of the $\ker(\lambda)$, 
$\lambda\in\Lambda_X$.

By Proposition \ref{tgt2} (ii), $T^1(X)^T=0$ for any
multiplicity-free $T$-variety $X$ (in fact, it is easy 
to show that every family of such varieties is locally 
trivial, see e.g. \cite[Lemma 7.4]{Alexeev_BrionI}; 
we refer to \cite{Altmann94} and \cite{Altmann97} for a study 
of arbitrary deformations of affine toric varieties). 
So we obtain the following result, which generalizes 
\cite[Theorem 1.2]{Peeva_Stillman02}.

\begin{corollary}\label{tor}
Let $X$ be an orbit closure in a finite-dimensional 
$T$-module $V$. If $X$ is non-degenerate, then the orbit 
$\GL(V)^T \cdot X$ is open in $\Hilb^T_{\cS}(V)$, where $\cS$ 
denotes the submonoid of $\Lambda$ generated by the opposites 
of the weights of $V$.
\end{corollary}

If, in addition, $X$ spans $V$, then all weight subspaces
have dimension $1$, so that $\GL(V)^T$ is a maximal torus of 
$\GL(V)$. In that case, $\Hilb^T_{\cS}(V)$ is the toric Hilbert 
scheme of \cite{Peeva_Stillman02}, and the closure of 
$\GL(V)^T \cdot X$ is its main component.

We now extend this construction to multiplicity-free 
$G$-varieties. Consider a finite-dimensional $G$-module
$V$, and denote by $F$ the set of weights of $(V^*)^U$.
Let $\cS$ the submonoid of $\Lambda^+$ generated by $F$,
and let $\M_{\cS}=\M_Y$ where $Y=\Spec k[\cS]$; 
then $Y$ is a non-degenerate, multiplicity-free subvariety 
of the $T$-module $V_U$.

\begin{corollary}\label{sph}
With this notation, the non-degenerate
multiplicity-free subvarieties $X\subseteq V$, having $\cS$ 
as their weight monoid, are parameterized by an open subscheme 
of $\Hilb^G_{\cS}(V)$, stable under $\GL(V)^G$. 
Moreover, this open subscheme $\Hilb^G_F$ is isomorphic to
$(\GL(V)^G\times \M_{\cS})/T$, where $T$ acts on $\GL(V)^G$ 
via $t\mapsto (\lambda(t))_{\lambda\in F}$, and on $\M_{\cS}$
via its action on $k[\cS]$.
\end{corollary}

\begin{proof}
Let $X\subseteq V$ be a multiplicity-free subvariety. If $X$ 
is non-degenerate, then so is its image in $V_U$.  
And if, in addition, the monoid $\Lambda^+_X=\Lambda^+_{X//U}$ 
is generated by $F$, then the map $X//U\rightarrow V_U$ 
is a closed immersion; its image is contained in the orbit of 
$Y$ under $\GL(V_U)^T\simeq \GL(V)^G$. Thus, $X$ is a closed 
point of the preimage in $\Hilb^G_{\cS}(V)_0$ of this orbit. 

Conversely, any closed point $X$ of $f^{-1}(\GL(V)^G\cdot Y)$
is a non-degenerate, multiplicity-free variety, with weight 
monoid $\cS$ (since $X//U\simeq Y$). So 
$\Hilb^G_F = f^{-1}(\GL(V)^G\cdot Y)$; then its 
structure follows from Theorem \ref{rep3}.
\end{proof}

A special $\GL(V)^G$-orbit in $\Hilb_F^G$ may be constructed
as follows. Choose $B$-eigenvectors $x_{\lambda}$ in all 
isotypical components $V_{(\lambda)}$, and let 
$x=\sum_{\lambda\in F} x_{\lambda}$. Then
$X_0=\overline{G\cdot x}$ is a non-degenerate $G$-subvariety
of $V$; one checks that it is also multiplicity-free with weight 
monoid $\cS$. 

Clearly, the $\GL(V)^G$-orbit of $X_0$ in $\Hilb^G_F$ is obtained 
by varying the choices of the $x_{\lambda}$, and the corresponding
$T$-orbit in $M_{\cS}$ is just a fixed point that we still
denote by $X_0$. 
We will see in Theorem \ref{aff} that the scheme $M_{\cS}$ is affine,
with $X_0$ as its unique closed orbit of $T$; as a consequence, 
$\Hilb_F^G$ is affine as well, and $\GL(V)^G\cdot X_0$ is its 
unique closed $\GL(V)^G$-orbit.

In the case where $G$ is a torus, $\Hilb^G_F$ is a unique orbit 
of $\GL(V)^G$. But this may fail in the general case,
as shown by the following examples.

\medskip

\noindent
{\sl Example 1}. Let $G=\SL(2)$. The simple $G$-modules are 
the symmetric powers of the defining module $k^2$; they are
indexed by the non-negative integers. We determine 
$\Hilb^G_F$, where $F=\{n\}$; then 
$\Hilb_F^G \simeq \M_{\bN n}$ by Corollary \ref{sph}.

The space $V=V(n)$ has basis the monomials 
$x^n,x^{n-1}y,\ldots,y^n$, where $x^n$ is the highest weight
vector. So $X_0=\overline{G\cdot x^n}$ is a normal surface, 
the cone over the rational normal curve in $\bP^n$, and 
$X_0 - G\cdot x^n$ is just the origin. Thus, we may apply
Proposition \ref{tgt2}. Note that
$G_{x^n} = U \mu_n$, where $\mu_n$ denotes the group of $n$-th
roots of unity in $T\simeq \bG_m$; it acts on $V(n)$ via 
$t\cdot(x,y)=(tx,t^{-1}y)$. Hence $\fg\cdot x^n$ is 
spanned by $x^n,x^{n-1}y$, and $(V(n)/\fg\cdot x^n)^U$ is spanned
by the image of $x^{n-2}y^2$ if $n\ge 2$, whereas 
$(V(1)/\fg\cdot x)^U$ is zero. But $x^{n-2}y^2$ is fixed by $\mu_n$
if and only if $n$ divides $4$. So $T_{X_0}\Hilb^G_F$ is a line
if $n=2$ or $n=4$, and vanishes otherwise.

It follows that $\Hilb^G_F$ is the affine line if $n=2$ or $n=4$,
and is a point otherwise. The corresponding families are the
orbit closures of $x^2 + t y^2$ in $V(2)$, and of $x^4 + tx^2y^2$
in $V(4)$, where $t\in k$; for $t\ne 0$, they are isomorphic to 
$G/T$, resp.$G/N_G(T)$, where $N_G(T)$ denotes the normalizer of 
$T$.

These results also follow from work of Pinkham 
\cite[Chapter 8]{Pinkham74}
describing all (possibly non-invariant) deformations of $X_0$.

\medskip

\noindent
{\sl Example 2.} Let $G=\SL(4)$ and 
$F=\{\omega_1,\omega_2,\omega_3\}$, where each $\omega_i$ 
is the highest weight of the simple $G$-module $\wedge^i k^4$.
Then $\cS=\Lambda^+$, and $\Hilb^G_F=\M_{\Lambda^+}$. We have
$$
V= k^4 \times \wedge^2 k^4 \times \wedge^3 k^4.
$$
One checks that the multiplicity-free subvarieties of $V$
that span it and have $\Lambda^+$ as their weight monoid
are exactly the orbit closures $\overline{G\cdot x}$, where
$x$ is one of the following points:
$$
x_0 = (e_1, e_1 \wedge e_2, e_1 \wedge e_2 \wedge e_3),~
x_1(t) = (e_1, e_2 \wedge e_3, t e_1 \wedge e_2 \wedge e_3),~
$$
$$
x_2(t) = (e_1, e_1 \wedge e_4, t e_1 \wedge e_2 \wedge e_3),
x_{12}(t,u) = (e_1, t(e_1 \wedge e_4 + e_2 \wedge e_3), 
u e_1 \wedge e_2 \wedge e_3).
$$ 
Here $(e_1,e_2,e_3,e_4)$ is an arbitrary basis of $k^4$,
and $t,u$ are non-zero scalars.

Fixing this basis, we obtain representatives of the $T$-orbits
in $M_{\Lambda^+}$; the dimensions of these orbits are 
$0,1,1,2$. In particular, $X_0=\overline{G\cdot x_0}$ is the
$T$-fixed point. The tangent space 
$T_{X_0}M_{\Lambda^+} = T_{X_0} \Hilb^G_{\cS}(V)$ 
may be determined by using Proposition \ref{tgt2}: indeed,
one may check that $X_0$ is a normal variety, with boundary
$X_0 - G\cdot x_0$ of codimension $2$. Moreover, $G_{x_0}$ is 
a maximal unipotent subgroup of $G$, and the images in 
$V/\fg\cdot x_0$ of 
$$
(0,e_1 \wedge e_4,0),~(0,e_2\wedge e_3,0)
$$
are a basis of $(V/\fg\cdot x_0)^{G_{x_0}}$. Thus, 
$T_{X_0} M_{\Lambda^+}$ has dimension $2$.

It follows that $M_{\Lambda^+}$ is isomorphic to an affine 
plane where $T$ acts linearly with weights 
$\alpha_1+\alpha_2,\alpha_2+\alpha_3$ where 
$\alpha_1,\alpha_2,\alpha_3$ are the simple roots. The 
$T$-orbit closures of $X_1(1)$, $X_2(1)$ are coordinate 
lines, and the orbit closure of $X_{12}(1,1)$ is the complement 
of their union.


\section{The action of the adjoint torus}

\label{sec: action}


\subsection{Definition of the action}

\label{subsec: definition}

Given $\lambda\in\Lambda$, we denote 
$e^{\lambda}\in k[T]$ the corresponding regular function on $T$. 
Then
$$
k[T]=\bigoplus_{\lambda\in \Lambda} k e^{\lambda}\supseteq 
\bigoplus_{\lambda\in\bZ\Pi} k e^{\lambda} = k[T_{\ad}]\supseteq
\bigoplus_{\lambda\in\bN\Pi} k e^{\lambda} = k[\bA^{\Pi}],
$$
where $\bA^{\Pi}$ is the affine space where $T$ acts linearly with
weights being the simple roots. The latter inclusion yields an open
immersion 
$$
T_{\ad}\hookrightarrow \bA^{\Pi},~t\mapsto (\alpha(t))_{\alpha\in\Pi}.
$$

We may identify the center $Z(G)$ with a central subgroup of $G\times T$,
via $z\mapsto(z,z^{-1})$. The quotient $\wG=(G\times T)/Z(G)$ is a
connected reductive group, with maximal torus 
$\wT =(T\times T)/Z(G)$. In fact, $\wT$ is isomorphic to 
$T\times T_{\ad}$ via $(t_1,t_2)Z(G)\mapsto (t_1t_2,t_2Z(G))$; this
defines an injective homomorphism 
$$
T_{\ad}\hookrightarrow \wT,~tZ(G)\mapsto (t^{-1},t)Z(G).
$$ 
On the other hand, we have an injective homomorphism 
$$
G\hookrightarrow \wG,~g\mapsto (g,1)Z(G).
$$ 
Together, these yield an isomorphism of $\wG$ with the semi-direct
product of $G$ with $T_{\ad}$, acting on $G$ by conjugation.

We now extend this construction to families of affine
$G$-schemes. Let $\pi:\cX\rightarrow S$ be such a family, then $Z(G)$
acts on $\cX\times T$ by $z\cdot (x,t)=(z\cdot x,z^{-1}t)$.  This
action is free, and $\pi$ is affine and $Z(G)$-invariant; thus, the 
quotient
\begin{eqnarray*}
\wcX=(\cX\times T)/Z(G)
\end{eqnarray*} 
is a scheme. The map 
$\cX\times G\times T\to \cX\times T,~(x,g,t)\mapsto (g\cdot x,t)$
descends to a morphism $\cX\times\wG \to \wcX$ which is
$G$-invariant, and the induced morphism
$(\cX\times \wG)/G\to \wcX$ is an isomorphism. Thus, $\wcX$ is
equipped with an action of $\wG$, and hence of $G$.
The morphism $\pi\times\id:\cX\times T\rightarrow S\times T$ is
affine, of finite type, and $Z(G)$-invariant. Hence it descends to an
affine morphism of finite type
$$
\wpi:\wcX\rightarrow (S\times T)/Z(G) = S\times T_{\ad}
$$
which is $G$-invariant and $T_{\ad}$-equivariant. In other words,
$\wpi$ is a family of affine $G$-schemes, with a compatible action of
$T_{\ad}$. Putting $\wcR=\wpi_*\cO_{\wcX}$, with isotypical components
the $\cR_{(\lambda)}$, we obtain
$$
\wcR=(\cR\otimes_k k[T])^{Z(G)}
=\bigoplus_{\lambda\in\Lambda^+,\mu\in\Lambda}
(\cR_{(\lambda)} e^{\mu})^{Z(G)}
=\bigoplus_{\lambda\in\Lambda^+,\mu\in\Lambda,\lambda-\mu\in\bZ\Pi}
\cR_{(\lambda)} e^{\mu}.
$$
This yields an isomorphism of $\cO_S[T_{\ad}]$-$G$-modules
\begin{eqnarray*}
\wcR = \bigoplus_{\lambda\in\Lambda^+}
\cR_{(\lambda)} e^{\lambda} \otimes_k k[T_{\ad}].
\end{eqnarray*}
Moreover, $\wcX//U = (\cX//U\times T)/Z(G)$ is mapped isomorphically
to $\cX//U\times T_{\ad}$, via $(x,t)\mapsto (t\cdot x,tZ(G))$.
In other words, the family 
$\wpi//U:\wcX//U\rightarrow S\times T_{\ad}$
is equipped with an isomorphism to the pull-back of 
$\pi//U:\cX//U\rightarrow S$ via the projection 
$S\times T_{\ad}\rightarrow S$.

If the family $\pi$ is flat, then so is $\wpi$, and its fiber at
any $(s,t)\in S\times T_{\ad}$ is $G$-equivariantly isomorphic to the
fiber of $\pi$ at $s$. Moreover, if $\pi$ is multiplicity-finite,
then so is $\wpi$, with the same Hilbert function $h$. 

In the case where $\cX=V(\lambda)$ where $G$ acts linearly (and
$\pi$ is constant), then $Z(G)$ acts by the restriction of the
character $\lambda$. Then, for any $\mu\in\Lambda$ such that 
$\mu\vert_{Z(G)}=\lambda\vert_{Z(G)}$ (that is,
$\mu-\lambda\in\bZ\Pi$), we obtain an isomorphism
$$
\begin{matrix}
\widetilde{V(\lambda)} =(V(\lambda)\times T)/Z(G) 
& \to & V(\lambda)\times T_{\ad} \cr
(v,t)Z(G) & \mapsto & (\mu(t)v,tZ(G)),\cr
\end{matrix}
$$
which is $\wG$-equivariant for the action of $G$ on 
$V(\lambda)\times T_{\ad}$ via $g\cdot (v,s) = (g\cdot v,s)$, and the
action of $T_{\ad}$ via $t\cdot (v,s)= (\mu(t)t^{-1}v, ts)$. We choose
$\mu=w_o\lambda$, so that $T_{\ad}$ acts on $V(\lambda)$ by 
\begin{eqnarray*}
t\cdot v = (w_o\lambda)(t) t^{-1} v.
\end{eqnarray*}
Then $V(\lambda)^{U^-}$ is fixed pointwise, and 
$q:V(\lambda)\to V(\lambda)_U$ is invariant. 

More generally, any finite-dimensional $G$-module $V$ becomes a
$\wG$-module by letting $T_{\ad}$ act on every submodule $V(\lambda)$
as above, so that $V^{U^-}$ and $q:V\to V_U$ are $T_{\ad}$-invariant. 
This yields an $\wG$-equivariant isomorphism 
${\widetilde V} \to V\times T_{\ad}$. Thus, if $\cX$ is a closed
$G$-stable subscheme of $V\times S$, then $\wcX$ is a closed
subscheme of 
$$
(V\times S\times T))/Z(G)\simeq V\times T_{\ad}\times S,
$$
stable under $\wG$, and hence under $G$. So, given 
$(\pi:\cX\rightarrow S)$ in $\cHilb^G_h(V)(S)$, we obtain 
$(\wpi:\wcX\rightarrow T_{\ad}\times S)$ in
$\cHilb^G_h(V)(T_{\ad}\times S)$. 
Applying this to the universal family 
$\pi: \Univ^G_h(V) \rightarrow \Hilb^G_h(V)$,
we obtain an element of $\cHilb^G_h(V)(T_{\ad} \times \Hilb^G_h(V))$,
that is, a morphism of schemes
\begin{eqnarray*}
a:T_{\ad}\times \Hilb^G_h(V) \rightarrow \Hilb^G_h(V).
\end{eqnarray*}

\begin{proposition}\label{act}
(i) With the preceding notation, $a$ is an action of $T_{\ad}$ on
$\Hilb^G_h(V)$. The open subscheme $\Hilb^G_h(V)_0$ is stable under
$T_{\ad}$, and the morphism 
$f:\Hilb^G_h(V)_0 \to \Hilb^T_h(V_U)$ is invariant.

(ii) The restriction 
$T_{\ad}\times \Hilb^G_h(V)_0\rightarrow \Hilb^G_h(V)_0$ extends
uniquely to a morphism 
$$
\bA^{\Pi}\times \Hilb^G_h(V)_0\rightarrow \Hilb^G_h(V)_0.
$$
Hence this morphism is an action of the multiplicative monoid
$\bA^{\Pi}$, leaving $f$ invariant.
\end{proposition}

\begin{proof} (i) follows from the preceding discussion; in fact, $a$ 
comes from the action of $T_{\ad}$ on $V$, normalizing the
actions of $G$, $U$, and fixing $q$.

(ii) Let ($\pi:\cX\rightarrow S$) in $\cHilb^G_h(V)_0(S)$. For $\cR$
and $\wcR$ as above, consider the sheaf 
$$
\widehat \cR = 
\bigoplus_{\lambda\in\Lambda^+,\mu\in\Lambda,\lambda-\mu\in\bN\Pi}
\cR_{(\lambda)} e^{\mu}
\subseteq  
\bigoplus_{\lambda\in\Lambda^+,\mu\in\Lambda,\lambda-\mu\in\bZ\Pi}
\cR_{(\lambda)} e^{\mu}=\wcR.
$$
Since any simple $G$-submodule $V(\nu)$ of a tensor product
$V(\lambda)\otimes_k V(\mu)$ satisfies 
$\lambda+\mu-\nu\in\bN\Pi$, it follows that $\widehat \cR$ is
a sheaf of $\cO_S[\bA^{\Pi}]$-$G$-subalgebras of $\wcR$. Moreover,
$$
\widehat \cR = 
\bigoplus_{\lambda\in\Lambda^+}
\cR_{(\lambda)} e^{\lambda}\otimes_k k[\bA^{\Pi}].
$$
Thus, $\widehat\cR$ is locally free over $\cO_S[\bA^{\Pi}]$, and
${\widehat\cR}^U$ is isomorphic to $\cR^U[\bA^{\Pi}]$. This defines a 
flat family of affine $G$-schemes 
$$
\widehat\pi:\widehat\cX\rightarrow S\times \bA^{\Pi}
$$
extending $\wpi:\wcX\rightarrow S\times T_{\ad}$, such that
$\widehat\pi//U$ is equipped with an isomorphism to the pull-back of
$\pi//U$. By assumption, $\cX\subseteq V\times S$, and the composition
$\cX//U\rightarrow V//U\times S\rightarrow V_U\times S$ is a closed
$T$-equivariant immersion. This yields a closed $T$-equivariant
immersion 
${\widehat\cX}//U\simeq \cX//U\times\bA^{\Pi}\hookrightarrow
V_U\times S\times\bA^{\Pi}$, 
and hence a closed $G$-equivariant immersion 
${\widehat\cX} \hookrightarrow V\times S\times \bA^{\Pi}$
extending $\wcX \hookrightarrow V\times S\times T_{\ad}$.
\end{proof}

The above Proposition defines an action of $T_{\ad}$ on $\M_Y$ 
(where $Y$ is a multiplicity-finite $T$-scheme), that extends to
$\bA^{\Pi}$. Here is a more direct construction of this action.

\begin{lemma}\label{com}
Consider the action of $T$ on $\M_Y$ via the homomorphism 
$T \rightarrow \Aut^T(Y)$. Then the induced action of $Z(G)$ is 
trivial, and the corresponding action of $T/Z(G)$ is the same
as the preceding $T_{\ad}$-action.
\end{lemma}

\begin{proof}
Let $(\pi:\cX\rightarrow S,\varphi:\cX//U\rightarrow Y\times S)$
be a family of type $Y$ and let $\gamma \in \Hom(S,Z(G))$ act on
$Y\times S$. Then $\gamma$ defines $\psi\in\Aut^G_S(\cX)$, such
that $\gamma \circ \varphi = \varphi \circ (\psi//U)$. Thus,
the family $(\pi,\gamma\circ\varphi)$ is equivalent to 
$(\pi,\varphi)$, so that $Z(G)$ acts trivially on $\M_Y$.

Choose a closed $T$-equivariant immersion $Y\hookrightarrow V_U$, 
where $V$ is a finite-dimensional $G$-module. By construction,
the action of $T_{\ad}$ on $\Hilb^G_h(V)$ arises from the 
$T$-action on $V^*$, where each isotypical component 
$V^*_{(\lambda)}$ is the eigenspace of weight $\lambda$. This 
action stabilizes $(V^*)^U=(V_U)^*$, and it restricts to the 
natural $T$-action on that space, hence on $V_U$ and on $Y$. 
This completes the proof.
\end{proof}

\subsection{Horospherical schemes as fixed points}

\label{subsec: horospherical}
We will construct a $T_{\ad}$-invariant section of the morphism 
$f:\Hilb^G_h(V)_0\to \Hilb^T_h(V_U)$. For this, we need some
preliminaries on representation theory.

Let $E$ be a rational $B$-module. Put 
$$
\Ind_B^G(E)=\Hom^B(G,E)=(k[G]\otimes_k E)^B,
$$ 
where $B$ acts by right multiplication on $k[G]$, and 
simultaneously on $E$. Then the action of $G$ on $k[G]$ by left 
multiplication equips $\Ind_B^G(E)$ with the structure of a 
rational $G$-module, the \emph{induced module} of $E$. 
The evaluation at the identity $e\in G$ yields a 
$B$-equivariant morphism 
$$
\varepsilon: \Ind_B^G(E) \rightarrow E,
$$
which is universal for $B$-equivariant morphisms from rational
$G$-modules to $E$. If, in addition, $E$ is a $B$-algebra, then so
is $\Ind_B^G(E)$, and $\varepsilon$ is an algebra homomorphism.

In the case where $E=k_{\lambda}$, the one-dimensional $B$-module
with weight $\lambda\in\Lambda$, then $\Ind_B^G(E) = V(-\lambda)^*$ 
if $\lambda\in -\Lambda^+$; otherwise, $\Ind_B^G(E)=0$. As a
consequence, if $E$ is finite-dimensional, then so is $\Ind_B^G(E)$.

Recall that any rational $G$-module $V$ defines the space $V_U$ of its
co-invariants, a rational $B$-module with trivial action of $U$ and
weights in $-\Lambda^+$. The $B$-equivariant map $q:V\rightarrow V_U$
factors uniquely through a $G$-equivariant map
$$
V\rightarrow \Ind_B^G(V_U),~v\mapsto (g\mapsto q(g\cdot v))
$$
which is in fact an isomorphism. The assignments
$E\mapsto\Ind_B^G(E)$ and $V\mapsto V_U$ yield an equivalence of the
category of rational $G$-modules, with the category of rational
$B$-modules with trivial $U$-action and weights in $-\Lambda^+$.

We define the (restricted) dual $V^*$ of a rational module $V$ as
the largest rational submodule of the space of linear forms on
$V$. Then we define the \emph{co-induced module} of a rational
$B$-module $E$: 
\begin{eqnarray*}
\Coind_B^G(E)=(\Ind_B^G(E^*))^*.
\end{eqnarray*}
This is a rational $G$-module, equipped with a
$B$-equivariant map 
$$
\iota:E\rightarrow \Coind_B^G(E)
$$ 
which is universal for maps from $E$ to rational $G$-modules. Note
that $\Coind_B^G(k_{\lambda})$ equals $V(\lambda)$ if $\lambda\in
\Lambda^+$, and vanishes otherwise. Moreover, $\Coind_B^G(E)$ is 
finite-dimensional if $E$ is. Note also that 
$$
(V^U)^*= (V^*)_U.
$$ 
So the assignments $V\rightarrow V^U$ and 
$E\rightarrow \Coind_B^G(E)$ yield an equivalence of the category
of rational $G$-modules, with the category of rational $T$-modules
with weights in $\Lambda^+$. 

Another relation between induced and co-induced modules arises from
the isomorphism $V^{U^-}\rightarrow V_U$, where $V$ is any rational
$G$-module. This implies an isomorphism
$$
\Coind_{B^-}^G(E)\simeq \Ind_B^G(E),
$$
where $E$ is any rational $T$-module (regarded as a trivial
$U^-$-module on the left-hand side, and as a trivial $U$-module on
the right-hand side). If $E$ is a $T$-algebra, this yields a
$G$-algebra structure on $\Coind_{B^-}^G(E)$.

Next we extend these constructions to families. Let
$\pi:\cY\rightarrow S$ be a family of affine $T$-schemes. Put 
$$
\cA=\pi_*\cO_{\cY} \text{ and } \cR=\Coind_B^G(\cA),
$$
where $B$ acts on $\cA$ via its quotient $T$; then the canonical map
$\iota:\cA\rightarrow \cR$ is $U$-invariant. Note that $\cR$ only
depends on $\bigoplus_{\lambda\in\Lambda^+} \cA_{\lambda}$, a sheaf of
finitely generated $\cO_S$-$T$-subalgebras of $\cA$. Thus, we may and
will assume that all weights of $\cA$ lie in $\Lambda^+$. We have
$$
\cR\simeq \Ind_{B^-}^G(\cA)= (k[G]^{U^-}\otimes_k \cA)^T,
$$
where $T$ acts diagonally on $k[G]^{U^-}\otimes_k \cA$. It follows
that $\cR$ is a sheaf of finitely generated
$\cO_S$-$G$-algebras; $\iota:\cA\rightarrow\cR$ is an injective
homomorphism of $\cO_S$-$B$-algebras identifying $\cA$ with $\cR^U$;
and $\varepsilon:\cR\to\cA$ is a surjective homomorphism of
$\cO_S$-$B^-$-algebras identifying $\cA$ with $\cR_{U^-}$.

\begin{definition}
A family of affine $G$-schemes $\pi:\cX \to S$ is \emph{horospherical}
if $\pi_*\cO_{\cX}=\Coind_B^G(\cA)$ for some sheaf $\cA$ of finitely
generated $\cO_S$-$T$-algebras with weights in $\Lambda^+$. 
Denoting $\cY=\Spec_S(\cA)$, a family of affine $T$-schemes, we say
that $\cX$ is \emph{induced from $\cY$} and write $\cX=\Ind_B^G(\cY)$.
\end{definition}

Then $\cX$ is equipped with an isomorphism $\cX//U\to \cY$. Further,
the map
$\Aut^G_S(\cX) \to \Aut^T_S(\cX//U)=\Aut^T_S(\cY)$ 
is an isomorphism. In particular, for any multiplicity-finite 
$T$-scheme $Y$, we may regard $\Ind_B^G(Y)$ as a closed point 
of $\M_Y$, fixed by $\Aut^T(Y)$. 

In \cite{Knop90}, a $G$-variety $X$ is said to be horospherical 
if $X=G\cdot X^U$. We now show that both definitions agree.

\begin{lemma}\label{hor}
The following conditions are equivalent, for a family of affine
$G$-schemes $\pi:\cX\to S$ with $U^-$-fixed point subscheme
$\cX^{U^-}$: 

\noindent
(i) $\cX$ is horospherical.

\noindent
(ii) $G\cdot \cX^{U^-}=\cX$ (as schemes).

\noindent
(iii) The product in $\cR=\pi_*\cO_{\cX}$ of any two isotypical
components $\cR_{(\lambda)}$, $\cR_{(\mu)}$ is contained in
$\cR_{(\lambda+\mu)}$. 

Then the map $\cX^{U^-} \to \cX//U$ is an isomorphism.
\end{lemma}

\begin{proof} We may assume that $\cX$ is affine and we put
$R=\Gamma(\cX,\cO_{\cX})$. Note that 
$G\cdot \cX^{U^-} = G\cdot \cX^U$ (since $U^-$ and $U$ are conjugate
in $G$). Further, $\cX^U$ is a closed $B$-stable subscheme of $\cX$, 
and hence $G\cdot \cX^U$ is also closed in $\cX$ (since $G/B$ is
complete).

(i)$\Rightarrow$(ii) We have $R\simeq \Ind_{B^-}^G(R_{U^-})$, so that
the quotient map $q^-:R\rightarrow R_{U^-}$ is an algebra homomorphism:
its kernel $\ker(q^-)$ is an ideal of $R$. But the ideal of the
fixed point subscheme $X^{U^-}$ is generated by the $g\cdot f - f$ for
$g\in U^-$ and $f\in R$; hence this ideal equals $\ker(q^-)$. Now
$\ker(q^-)$ contains no simple $G$-submodule of $R$, so that 
$\bigcap_{g\in G} g\cdot \ker(q^-)=\{0\}$. It follows that 
$G\cdot X^{U^-}=X$.

(ii)$\Rightarrow$(iii) The group $G\times B^-$ acts on 
$G\times X^{U^-}$ via $(g_1,\gamma)\cdot(g_2,y) =
(g_1 g_2 \gamma^{-1},\gamma\cdot y)$, and the surjective 
morphism $G\times X^{U^-}\rightarrow X$ is $G$-equivariant 
and $B^-$-invariant. Thus, it yields an injective homomorphism 
of $G$-algebras
$$
R \hookrightarrow (k[G]^{U^-}\otimes_k A)^T,
$$
where $A=\Gamma(X^{U^-},\cO_{X^{U^-}})$.
Since each weight space $k[G]^{U^-}_{\lambda}$ is a simple $G$-module 
(with highest weight $-\lambda$), it follows that any simple
$G$-submodule of $R$ can be written as 
$k[G]^{U^-}_{-\lambda}\otimes f$, for some $\lambda\in \Lambda^+$ and
some $f\in A_{\lambda}$. Now (iii) follows from the fact that
$$
k[G]^{U^-}_{-\lambda} \; k[G]^{U^-}_{-\mu} =
k[G]^{U^-}_{-\lambda-\mu}
$$
for all $\lambda$, $\mu$ in $\Lambda^+$.

(iii)$\Rightarrow$(i) It suffices to show that the kernel $\ker(q^-)$ 
of the quotient map $q^-:R\rightarrow R_{U^-}$ is an ideal of $R$,
i.e., is stable under multiplication by any simple $G$-submodule
$V(\mu)\subseteq R$.
For this, note that the intersection of $\ker(q^-)$ with any simple
$G$-submodule $V(\lambda)\subseteq R$ is the kernel of
$V(\lambda)\mapsto V(\lambda)_{U^-}$, that is, the sum of all weight
subspaces $V(\lambda)_{\chi}$ where $\chi\neq \lambda$; then
$\chi<\lambda$. Thus, $\chi + \eta < \lambda + \mu$, for any weight
$\eta$ of $V(\mu)$. Since the product $V(\lambda)V(\mu)$ is
either $0$ or $V(\lambda+\mu)$, this implies our assertion.
\end{proof}

Next we relate the induction $\cY\mapsto\Ind_B^G(\cY)$ to invariant
Hilbert schemes.

\begin{proposition}\label{sec}
(i) The induction of families of affine $T$-schemes yields a closed
immersion 
\begin{eqnarray*}
s:\Hilb^T_h(V_U) \rightarrow \Hilb^G_h(V)_0
\end{eqnarray*}
which is a section of $f:\Hilb^G_h(V)_0 \rightarrow \Hilb^T_h(V_U)$.

\noindent
(ii) The image of $s$ is contained in the fixed point
subscheme of $T_{\ad}$ in $\Hilb^G_h(V)_0$. Thus, $s$ is invariant
under $\bA^{\Pi}$. 

\noindent
(iii) For any closed point $Y$ of $\Hilb^T_h(V_U)$, the origin
$o\in\bA^{\Pi}$ acts on $\M_Y=f^{-1}(Y)$ as the constant map 
to the closed point $\Ind_B^G(Y)$.
\end{proposition}

\begin{proof}
(i) Let $\pi:\cX\to S$ be induced from a family 
$\rho:\cY\to S$ of affine $T$-schemes, such that all weights of
$\rho_*\cO_{\cY}$ are in $\Lambda^+$. Then $\pi$ is flat if and only
if $\rho=\pi//U$ is. If, in addition, $\rho$ is multiplicity-finite
with Hilbert function $h$, then $\pi$ is multiplicity-finite with 
Hilbert function $h\vert_{\Lambda^+}$. Finally, if $\cY$ is a closed
$T$-stable subscheme of $V_U\times S$, then the $B$-equivariant map
$(V_U)^*\to \cA$ yields a $G$-equivariant map
$$
V^*=\Coind_B^G((V_U)^*) \rightarrow \cR=\pi_*\cO_{\cX},
$$
and hence a homomorphism 
$\cO_S\otimes_k \Sym_k(V^*)\rightarrow \cR$. 
The latter is surjective, since the corresponding homomorphism 
$$
\cO_S\otimes_k \Sym_k(V^*)^U = \cO_S\otimes_k \Sym_k(V_U)^*
\rightarrow \cR^U=\cA
$$ 
is. So we obtain a closed $G$-equivariant immersion 
$\cX\hookrightarrow V\times S$
lifting the immersion $\cY\hookrightarrow V_U\times S$, and hence an
element of $\cHilb^G_h(V)_0(S)$ mapped to $\cY$ under 
$f:\cHilb^G_h(V)_0(S) \rightarrow \cHilb^T_h(V_U)(S)$. 
This constructs the section $s$, which is automatically a closed
immersion.

(ii) With the preceding notation, one has by Lemma \ref{hor}:
$$
\cX=G\cdot \cX^{U^-}\subseteq G\cdot V^{U^-}\times S,
$$
and $T_{\ad}$ acts on $V$ by normalizing the $G$-action and fixing 
$V^{U^-}$ pointwise. Thus, $\cX$ is invariant under $T_{\ad}$.

(iii) Let $\pi:\cX\to S$ be an $S$-point of $\M_Y$. By definition of
the action, $o$ maps this point to the pullback to 
$S\times \{o\}$ of the family 
$\widehat\pi : \widehat\cX \rightarrow S\times\bA^{\Pi}$
constructed in the proof of Proposition \ref{act}. 
Let $o\cdot\cX$ be the corresponding scheme, and $o\cdot \cR$ 
the corresponding sheaf of $\cO_S$-$G$-algebras. Since 
$$
{\widehat\pi}_*\cO_{\widehat\cX}=
\bigoplus_{\lambda\in\Lambda^+} \cR_{(\lambda)} e^{\lambda}\otimes_k
k[\bA^{\Pi}]
$$ 
as sheaves of $\cO_S[\bA^{\Pi}]$-$G$-modules, and 
$$
\cR_{(\lambda)} e^{\lambda}\cdot \cR_{(\mu)} e^{\mu}
\subseteq
\cR_{(\lambda+\mu)} e^{\lambda+\mu} +
\sum_{\gamma\in\bN\Pi,\gamma\ne 0}
\cR_{(\lambda+\mu-\gamma)} e^{\lambda+\mu-\gamma}
e^{\gamma},
$$
we see that $o\cdot \cR$ satisfies condition (iii) of Lemma \ref{hor}.
By that Lemma, $o\cdot\cX$ is horospherical. Now the
isomorphism $o\cdot\cX//U \rightarrow Y\times S$ implies that
$o\cdot\cX=\Ind_B^G(Y)\times S$.
\end{proof}


\subsection{Structure and tangent space of $\M_Y$}

\label{subsec: structure}

We begin with a general result on torus actions which is probably
known, but for which we could not find a reference. 

\begin{lemma}\label{lim}
Let $\bG_m^n\times X\rightarrow X,~(t,x)\mapsto t\cdot x$ be an
action of the $n$-dimensional torus $\bG_m^n$ on a scheme $X$, which
extends to a morphism $\bA^n\times X\rightarrow X$ under the natural
inclusion $\bG_m^n\subset\bA^n$; let $o$ be the origin of
$\bA^n$. Then the image of the morphism 
$o:X\rightarrow X,~x\mapsto o\cdot x$ is the fixed point subscheme
$X^{\bG_m^n}$, and the corresponding morphism  
$f:X\rightarrow X^{\bG_m^n}$ is a good quotient by $\bG_m^n$. 
\end{lemma}

\begin{proof}
In the case where $X=\Spec(R)$ is affine, the assumption that the
action extends to $\bA^n$ means that all weights of $\bG_m^n$ in
$R$ lie in $\bN^n$. Then the projection 
$R\rightarrow R^{\bG_m^n}=R_0$ is an algebra homomorphism; regarded as
a self-map of $R$, it is the comorphism of $o$. This implies easily
our statements.

In the general case, note that the image of $o$ is
contained in $X^{\bG_m^n}$. Moreover, the resulting morphism
$f:X\rightarrow X^{\bG_m^n}$ restricts to the identity on
$X^{\bG_m^n}$; in particular, $f$ is surjective. By the first step of
the proof, it suffices to show that $f$ is affine. For this, we may
assume that $X$ is a variety.

Let $\nu:\wX\rightarrow X$ be the normalization map. Then the map
$$
\bA^n\times\wX\rightarrow X,~(z,\xi)\mapsto z\cdot\nu(\xi)
$$
is dominant, and hence factors through a morphism
$\bA^n\times\wX\rightarrow \wX$. The latter is an action of the
multiplicative monoid $\bA^n$, since $\nu$ is birational; this yields
a morphism 
${\widetilde f}: \wX\rightarrow \wX^{\bG_m^n}$ lifting $f$. And since
$\nu$ is finite surjective, its fiber at any fixed point is stable 
under this action, and hence under $\bG_m^n$. Thus, $\nu$ restricts 
to a finite surjective morphism 
$\eta:\wX^{\bG_m^n}\rightarrow X^{\bG_m^n}$.
 
By a theorem of Chevalley \cite[Exercise III.4.2]{Hartshorne77}, 
an open subset of $X$ is affine if and only if its preimage under
$\nu$ is. Thus, it suffices to show that the composition 
$f\circ \nu$ is affine. Since 
$f\circ \nu = \eta \circ {\widetilde f}$,
we may further assume that $X$ is normal. 

By \cite{Sumihiro74}, $X$ is covered by open affine
$\bG_m^n$-stable subsets $X_i$. Then the $X_i^{\bG_m^n}$
form an open affine covering of $X^{\bG_m^n}$. Thus, to show that
$o$ is affine, it suffices to check that the preimage of each
$X_i^{\bG_m^n}$ is affine. Note that this preimage is contained in
$X_i$ (namely, $o\cdot x\in X_i$ implies that the closure of
$\bG_m^n\cdot x$ meets $X_i$, whence $x\in X_i$ since $X_i$ is open
and $\bG_m^n$-stable). Thus, replacing $X$ with 
$o^{-1}(X_i^{\bG_m^n})$, we may assume that: $X$ is contained in an
affine $\bG_m^n$-variety $V$ as a locally closed $\bG_m^n$-stable
subvariety, and $X^{\bG_m^n}$ is closed in $V$. We may further assume
that $V$ is a $\bG_m^n$-module. 

Let $V_0=V^{\bG_m^n}$ and $V_+=\{v\in V~\vert~ o\cdot v = 0\}$
(the sum of all ``positive'' $\bG_m^n$-eigenspaces). 
By our assumptions, $X^{\bG_m^n} = X\cap V_0$ is closed in $V_0$, and
$X$ is contained in $(X\cap V_0)\times V_+$ as a locally closed
$\bG_m^n$-stable subvariety. So the boundary $\oX-X$
is a closed $\bG_m^n$-stable subvariety of $(X\cap V_0)\times V_+$,
disjoint from $X\cap V_0$: this boundary must be empty. Hence $X$ is
closed in $V$ and, in particular, affine.
\end{proof}

Combining Propositions \ref{act} and \ref{sec} with Lemma
\ref{lim}, we obtain 

\begin{theorem}\label{aff}
The morphism $f:\Hilb^G_h(V)_0\rightarrow \Hilb^T_h(V_U)$
is a good quotient by the action of $T_{\ad}$, and the image of 
the section $s$ is the fixed point subscheme.

As a consequence, for any multiplicity-finite $T$-scheme $Y$, 
the scheme $\M_Y$ is affine and connected; its fixed point 
$\Ind_B^G(Y)$ is the unique closed orbit of $T_{\ad}$.
\end{theorem}

Next we describe the tangent space to $\M_Y$ at any closed point 
$X$. Since $X$ is equipped with an isomorphism $X//U\rightarrow Y$, 
we have restriction maps 
$$
r^0:\Der^G(X)\to \Der^T(Y) \text{ and } r^1:T^1(X)^G\to T^1(Y)^T.
$$
We may regard $r^0$ as the differential of the inclusion
$\Aut^G(X)\hookrightarrow \Aut^T(Y)$ (Lemma \ref{aut}).

\begin{proposition}\label{tgt3}
With the preceding notation, the restriction maps fit into an
exact sequence
\begin{eqnarray*}
0 \rightarrow  \Der^G(X) \rightarrow \Der^T(Y) \rightarrow
T_X \M_Y \rightarrow T^1(X)^G \rightarrow T^1(Y)^T \rightarrow 0.
\end{eqnarray*}
Thus, $\ker(r^1)$ identifies with the normal space at $X$ to
the orbit $\Aut^T(Y)\cdot X$ in $\M_Y$; in particular, if 
$X=\Ind_B^G(Y)$ then $\ker(r^1)=T_X \M_Y$.
\end{proposition}

\begin{proof}
With the notation of Theorem \ref{rep3}, consider the differential
$$
df_X:T_X\Hilb^G_h(V)\rightarrow T_Y\Hilb^T(V_U).
$$ 
It is surjective (since $f$ admits a section), with kernel 
$T_X f^{-1}(Y)=T_X \M_Y$. Moreover,
Proposition \ref{tgt1} yields a commutative diagram
$$
\begin{matrix}
0 \rightarrow & \Der^G(X) & \rightarrow & \Hom^G(X,V) & \rightarrow
& T_X\Hilb^G_h(V)_0 & \rightarrow & T^1(X)^G & \rightarrow 0 \cr
& \downarrow & & \downarrow & & \downarrow & & \downarrow &  \cr
0 \rightarrow & \Der^T(Y) & \rightarrow & \Hom^T(Y,V_U) & \rightarrow
& T_X\Hilb^T_h(V_U) & \rightarrow & T^1(Y)^T & \rightarrow 0\cr
\end{matrix}
$$
Here the map $\Hom^G(X,V)\rightarrow \Hom^T(Y,V_U)$ is the composition
\begin{eqnarray*}
\Hom^G(X,V)=\Hom^G(V^*,R)\rightarrow \Hom^T((V^*)^U,R^U)=\Hom^T(Y,V_U)
\end{eqnarray*}
(where $R=\Gamma(X,\cO_X)$, so that $R^U=\Gamma(Y,\cO_Y)$). Hence
$\Hom^G(X,V)$ is mapped isomorphically to $\Hom^T(Y,V_U)$; this implies
our exact sequence.
\end{proof}

If, in addition, $X$ is nonsingular, then $T^1(X)=0$, see e.g.
\cite[Example III.9.13.2]{Hartshorne77}. This yields

\begin{corollary}\label{rig}
Let $X$ be a nonsingular affine $G$-variety, and put $Y=X//U$.
Then the orbit $\Aut^T(Y)\cdot X$ is open in $\M_Y$.

As a consequence, there are only finitely many isomorphism
classes of nonsingular affine $G$-varieties of fixed type.
\end{corollary}

In particular, there are only finitely many isomorphism 
classes of nonsingular multiplicity-free $G$-varieties 
with a prescribed weight monoid. (This will be generalized
to possibly singular varieties in Corollary \ref{fin2} 
below.) In fact, a conjecture of Knop asserts that all
nonsingular multiplicity-free $G$-varieties are 
classified by their weight monoid. This conjecture has been 
established by Camus \cite{Camus01} for $G$ of type $A$, 
building on Luna's classification of spherical varieties 
for such $G$ \cite{Luna01}.


\subsection{Orbit closures}

\label{subsec: orbit}

Let $A$ be a multiplicity-finite $T$-algebra with weights in
$\Lambda^+$, and let $Y=\Spec(A)$. We will study the $T_{\ad}$-orbits
in the moduli space $\M_Y$; we begin with a more concrete description
of that space.

\begin{proposition}\label{str}
The functor $\cM_Y$ is isomorphic to the contravariant functor
(Schemes) $\rightarrow$ (Sets) assigning to $S$ the set of those 
$\cO_S$-$G$-algebra multiplication laws on 
$\cO_S\otimes_k\Coind_B^G(A)$, that extend the multiplication of $A$. 
\end{proposition}

\begin{proof}
Let $(\pi:\cX \rightarrow S,\varphi:\cX//U \rightarrow Y\times S)$
be a family of affine $G$-schemes of type $Y$. It yields an isomorphism
$$
\varphi^{\#}:\cO_S\otimes_k A \rightarrow \cR^U
$$
of $\cO_S$-$T$-algebras (where $\cR=\pi_*\cO_{\cX}$), which extends
uniquely to an isomorphism
$$
\Coind_B^G(\varphi^{\#}):
\cO_S\otimes_k \Coind_B^G(A) \rightarrow \cR =\Coind_B^G(\cR^U)
$$
of $\cO_S$-$G$-modules. Hence we obtain a structure of
$\cO_S$-$G$-algebra on $\cO_S \otimes_k \Coind_B^G(A)$, extending
the $T$-algebra structure of $A$. 

Clearly, any such structure arises from some family. Moreover, two
families $(\pi,\varphi)$ and $(\pi',\varphi')$ give rise to the same
structure if and only if there exists a $G$-equivariant isomorphism
$\psi:\cX \rightarrow \cX'$ over $S$ such that 
$\Coind_B^G(\varphi^{\#})\circ \psi^{\#}=\Coind_B^G(\varphi^{'\#})$,
that is, $\varphi=\varphi'\circ(\psi//U)$. In other words, these
families are equivalent in the sense of Definition \ref{typ}.
\end{proof}

Consider a $\cO_S$-$G$-algebra multiplication law
$$
m:\cR\otimes_{\cO_S}\cR\rightarrow \cR.
$$ 
Since $\cR=\bigoplus_{\lambda\in\Lambda^+} \cR_{(\lambda)}$, 
it follows that 
$$
m=\sum_{\lambda,\mu,\nu\in\Lambda^+}m_{\lambda,\mu}^{\nu}~,
$$
where every $m_{\lambda,\mu}^{\nu}$ lies in
$\Hom^G_{\cO_S}(\cR_{(\lambda)}\otimes_{\cO_S}\cR_{(\mu)},\cR_{(\nu)})$.
And since 
$\cR_{(\lambda)}\simeq \cO_S\otimes_k A_{\lambda}\otimes_k V(\lambda)$
as $\cO_S$-$G$-modules, we obtain
$$
m_{\lambda,\mu}^{\nu}\in \Gamma(S,\cO_S)\otimes_k
\Hom(A_{\lambda}\otimes_k A_{\mu},A_{\nu})\otimes_k
\Hom^G(V(\lambda)\otimes_k V(\mu),V(\nu))
$$
for all triples $(\lambda,\mu,\nu)$. 

As a consequence,
$m_{\lambda,\mu}^{\nu} =0$ unless $\nu\le \lambda+\mu$,
and $m_{0,\mu}^{\nu}=0$ unless $\nu = \mu$.  
Moreover, the condition that $m$ extends the multiplication 
of $\cO_S\otimes_k A$ is equivalent to
$$
m_{\lambda,\mu}^{\lambda+\mu} = 
1\otimes \mul_{\lambda,\mu} \otimes p_{\lambda,\mu}~,
$$
where $\mul_{\lambda,\mu}$ denotes the multiplication
$A_{\lambda}\otimes_k A_{\mu}\rightarrow A_{\lambda+\mu}$,
and 
$p_{\lambda,\mu} : V(\lambda)\otimes_k V(\mu) 
\rightarrow V(\lambda+\mu)$
is the unique morphism of $G$-modules sending 
$v_{\lambda}\otimes v_{\mu}$ to $v_{\lambda+\mu}$.
Under this condition, the identity element $1\in A$ is also
the identity for $m$. On the other hand, the commutativity
(resp. associativity) of $m$ translates into a family of 
linear (resp. quadratic) relations between the
$m_{\lambda,\mu}^{\nu}$'s.

Taking $S=\M_Y$, we see that any isotypical component 
$m_{\lambda,\mu}^{\nu}$ yields a morphism from $\M_Y$
to some affine space. For simplicity, we still
denote $m_{\lambda,\mu}^{\nu}$ this morphism, and $m$ 
the product of all the $m_{\lambda,\mu}^{\nu}$. Clearly, 
$m$ is universally injective. 

\begin{proposition}\label{coo}
Every isotypical component $m_{\lambda,\mu}^{\nu}$ is an 
eigenvector of $T_{\ad}$ of weight $\lambda+\mu-\nu$,
and these generate the algebra $\Gamma(\M_Y,\cO_{\M_Y})$.
\end{proposition}

\begin{proof}
By Lemma \ref{com}, the $T_{\ad}$-action on $\M_Y$ comes from
the $T$-action on 
$A=\bigoplus_{\lambda\in\Lambda^+} A_{\lambda}$. This implies 
the first assertion.

By that assertion or Lemma \ref{hor}, $m_{\lambda,\mu}^{\nu}$
vanishes at the $T_{\ad}$-fixed point $X_0=\Ind_B^G(Y)$
whenever $\lambda+\mu-\nu \ne 0$. Moreover, the tangent space 
$T_X \M_Y$ is spanned by the images of all such components,
since $m$ is universally injective. But the algebra 
$\Gamma(\M_Y,\cO_{\M_Y})$ is $\bN\Pi$-graded, with homogeneous
maximal ideal corresponding to $X_0$. So the
$m_{\lambda,\mu}^{\nu}$ generate this algebra, by the graded
Nakayama lemma.
\end{proof}

The above Proposition, together with Proposition \ref{str}, 
yields another, more direct proof of the second part of Theorem
\ref{aff}. It also motivates the following

\begin{definition}\label{roo}
Let $X=\Spec(R)$ be an affine $G$-scheme, and denote by
$m=\sum m_{\lambda,\mu}^{\nu} : R\otimes_k R\rightarrow R$
the multiplication. The \emph{root monoid} of $X$ is the 
submonoid $\Sigma_X\subseteq\Lambda$ generated by the 
$\lambda+\mu-\nu$, where $\lambda$, $\mu$, $\nu\in\Lambda^+$ 
and  $m_{\lambda,\mu}^{\nu}$ is nonzero.
\end{definition}

Note that $\Sigma_X$ is contained in $\bN\Pi$. Now Proposition
\ref{coo} easily implies

\begin{proposition}\label{mon}
For any affine $G$-scheme $X$, the root monoid $\Sigma_X$ is 
the weight monoid of the $T_{\ad}$-orbit closure of $X$
(regarded as a closed point of $\M_Y$). As a consequence,
the monoid $\Sigma_X$ is finitely generated.
\end{proposition}

If, in addition, $Y$ is a variety, then so is $X$. Let 
${\widetilde\Sigma}_X$ be the saturation of $\Sigma_X$, i.e.,
the intersection of the cone and the lattice generated
by that monoid (the corresponding multiplicity-free 
$T_{\ad}$-variety is the normalization of that associated 
with $\Sigma_X$). By \cite[Theorem 1.3]{Knop96}, the monoid 
${\widetilde\Sigma}_X$ is freely generated by a basis of a certain 
root system $\Phi_X$. Together with Propositions \ref{coo} 
and \ref{mon}, this implies

\begin{corollary}\label{nor} 
The normalization of every $T_{\ad}$-orbit in $\M_Y$ is 
isomorphic to a $T_{\ad}$-module, for any multiplicity-finite 
$T$-variety $Y$.
\end{corollary}


\section{Finiteness results for multiplicity-free varieties}

\label{sec: finiteness}

It is easy to see that only finitely many subgroups of 
a torus $T$ occur as isotropy groups of points of a given 
finite-dimensional $T$-module $V$; it follows that
$V$ contains only finitely many multiplicity-free 
subvarieties, up to the action of $\GL(V)^T$. In contrast,
there are examples of finite-dimensional $G$-modules 
where infinitely many pairwise non-isomorphic isotropy
groups occur, see \cite{Richardson68} 1.3 and 
\cite{Richardson72} 12.4.2. But finiteness still holds for
spherical isotropy groups:

\begin{theorem}\label{fin}
For any $G$-scheme $\cX$ of finite type, only finitely many
conjugacy classes of spherical subgroups of $G$ occur
as isotropy groups of points of $\cX$.
\end{theorem}

\begin{proof}
It proceeds through several reductions; we divide it into five 
steps.

\medskip

\noindent
\emph{Step 1.}
We reduce to the case of a nonsingular $G$-variety $\cX$ equipped 
with a smooth $G$-invariant morphism 
$$
\pi:\cX\rightarrow S
$$
such that the fiber $\cX_s$ at any closed point $s\in S$ is a 
unique $G$-orbit.

Indeeed, it suffices to prove the theorem for a $G$-variety $\cX$ 
(of finite type). By a classical result of Rosenlicht \cite{Ros63}, 
there exists a dense open subset $\cX_0\subseteq \cX$ admitting a
geometric quotient $\pi:\cX_0\rightarrow S$,
that is, $\pi$ is surjective, its fibers are the orbits, 
and the map $\cO_S\rightarrow (\pi_*\cO_{\cX_0})^G$ is 
an isomorphism. Shrinking $S$, we may further assume that 
it is nonsingular, and that $\pi$ is smooth. Thus $\cX_0$ 
is smooth as well. Now, using induction over the dimension, 
we may replace $\cX$ with $\cX_0$.

\medskip

\noindent
\emph{Step 2.}
After shrinking $S$ again, we may assume that $\cX$ admits a smooth
compactification $\ocX$, that is, a nonsingular $G$-variety 
containing $\cX$ as a dense open $G$-stable subset, such that 
$\pi$ extends to a smooth projective morphism 
$$
\overline{\pi} : \ocX \rightarrow S
$$
(then $\overline{\pi}$ is $G$-invariant).

Indeed, by \cite{Sumihiro75}, there exists a $G$-variety
$\ocX$ over $S$, containing $\cX$ as a dense open $G$-stable 
subset. Using equivariant desingularization 
\cite{Encinas_Villamayor00}, we may assume that $\ocX$ 
is nonsingular. Then, by generic smoothness, we may also assume
that $\overline{\pi}$ is smooth.

\medskip

\noindent
\emph{Step 3.}
After shrinking $S$ again, we may further assume that 
the compactification is ``regular'', that is, the boundary
$$
\partial \ocX = \ocX - \cX
$$
is a union of irreducible nonsingular divisors 
$\cY_1,\ldots,\cY_n$ with normal crossings, and each restriction
$$
\overline{\pi}_i : \cY_i \rightarrow S
$$
is smooth.

Indeed, we have an exact sequence
$$
0 \rightarrow \cE \rightarrow \cO_{\cX}\otimes_k\fg \rightarrow
\cT_{\cX/S} \rightarrow 0,
$$
where $\cT_{\cX/S}$ is the relative tangent sheaf and $\cE$ is 
a locally free sheaf on $\cX$ such that every fiber $\cE_x$ 
identifies with the isotropy Lie algebra $\fg_x$. This yields 
a morphism 
$$
\varphi:\cX \rightarrow \Grass(\fg),~x\mapsto \fg_x
$$
where $\Grass(\fg)$ denotes the Grassmannian of subspaces of 
$\fg$.

Replacing $\overline{\cX}$ with the closure in 
$\overline{\cX}\times\Grass(\fg)$ of the graph of $\varphi$, and 
resolving singularities, we may assume that $\varphi$ extends to 
$$
\overline{\varphi}:\overline{\cX} \rightarrow \Grass(\fg).
$$
So the relative tangent bundle $T_{\cX/S}$ extends to 
a vector bundle ${\widetilde T}_{\overline{\cX}/S}$ over $\ocX$, 
the pull-back of the tautological bundle of $\Grass(\fg)$. 
By \cite{Bien_Brion96}, it follows that every fiber $\oX=\ocX_s$ 
is a regular compactification of its open orbit $X=\cX_s$.
This implies easily our reduction.
 
\medskip

\noindent
\emph{Step 4.}
By \cite[Proposition 2.5]{Bien_Brion96}, the pull-back of 
${\widetilde T}_{\ocX/S}$ to any fiber $\oX$ is 
the logarithmic tangent bundle $T_{\oX}(-\log \partial\oX)$,
associated with the sheaf $\cT_{\oX}(-\log \partial\oX)$
of derivations of $\oX$ preserving the ideal sheaf of 
$\partial\oX$. The space of infinitesimal deformations of 
the pair $(\oX,\partial\oX)$ is the first cohomology group 
$H^1(\oX, \cT_{\oX}(-\log \partial\oX))$, see e.g.
\cite[Proposition 3.1]{Kawamata85}. And this group vanishes, 
as a special case of \cite[Theorem 4.1]{Knop94}. 
Now deformation theory should tell us that the equivariant
isomorphism class of $(\oX,\partial\oX)$ is independent of 
the fiber, which should complete the proof. But for lack of 
an appropriate reference, we will provide an alternative 
argument based on nested Hilbert schemes; for the latter, 
see \cite{Cheah98}.

Since the morphism $\overline{\pi}$ is projective, we may 
regard $\ocX$ as a $G$-stable subvariety of $\bP(V)\times S$, 
where $V$ is a finite-dimensional $G$-module. Replacing $V$
with a symmetric power, we may assume additionally that the 
restriction map
$$
V^*=H^0(\bP(V),\cO_{\bP(V)}(1))\rightarrow H^0 (X,\cO_X(1))
$$
is surjective for all fibers $X$.

Let $\Hilb$ be the nested Hilbert scheme that parameterizes 
those families $(Z,Z_1,\ldots,Z_n)$ of closed subschemes of $\bP(V)$ 
such that $Z_1\cup\cdots\cup Z_n\subseteq Z$, and the Hilbert 
polynomial of $Z$ (resp. $Z_1,\ldots,Z_n$) equals that of $\oX$ 
(resp. $Y_1,\ldots,Y_n$). Then we have a morphism
$\psi:S\rightarrow \Hilb$
such that $\ocX$ is the pull-back under $\psi$ of the universal 
family. The group $\GL(V)$, and hence $G$, acts on $\Hilb$, 
and $\psi$ is $G$-invariant. 

The differential of the $\GL(V)$-action yields a linear map
$$
f:\End(V) \rightarrow T_{(\oX,Y_1,\ldots,Y_n)}\Hilb
$$
which is $G$-equivariant. We will check in Step 5 below that 
$f$ is surjective. As a consequence, the restriction
$$
f^G:\End^G(V) \rightarrow (T_{(\oX,Y_1,\ldots,Y_n)} \Hilb)^G
$$
is surjective as well. Since the subscheme of $G$-invariants
$\Hilb^G$ is stable under the group $\GL(V)^G$ with Lie algebra 
$\End^G(V)$, and
$$
T_{(\oX,Y_1,\ldots,Y_n)}(\Hilb^G)\subseteq 
(T_{(\oX,Y_1,\ldots,Y_n)} \Hilb)^G,
$$
the orbit $\GL(V)^G \cdot (\oX,Y_1,\ldots,Y_n)$ is open in
$\Hilb^G$. It follows that all fibers of $\overline{\pi}$ 
are isomorphic in a neighborhood of $s$, as desired.

\medskip

\noindent
\emph{Step 5.}
To complete the proof, we deduce the surjectivity of $f$ from 
Knop's vanishing theorem stated in Step 4. For this, we recall
the description of the tangent space to the nested Hilbert scheme
\cite{Cheah98}.

Let $N_{\oX}$ (resp. $N_{Y_i}$) be the normal bundle to 
$\oX$ (resp. $Y_i$) in $\bP(V)$; then we have restriction maps 
$p_i: N_{Y_i}\rightarrow N_{\oX}\vert_{Y_i}$. 
Now $T_{(\oX,Y_1,\ldots,Y_n)}\Hilb$ equals
$$
\{(s,t_1,\ldots,t_n) \in H^0(N_{\oX}) \oplus
\bigoplus_{i=1}^n H^0(N_{Y_i})~\vert~ s\vert_{Y_i}=p_i(t_i)~ 
(i=1,\ldots,n)\}.
$$
We now construct a sheafified version of 
$T_{(\oX,Y_1,\ldots,Y_n)}\Hilb$, as follows. Denote by $\cN_{\oX}$ 
(resp. $\cN_{Y_i}$) the normal sheaf to $\oX$ (resp. $Y_i$) in
$\bP(V)$. We have exact sequences of sheaves on $Y_i$: 
$$
\CD
0 @>>> \cN_{Y_i/\oX} @>>> \cN_{Y_i} @>{p_i}>> 
\cN_{\oX}\vert_{Y_i} @>>> 0,
\endCD
$$
where $\cN_{Y_i/\oX}$ denotes the normal sheaf to $Y_i$ in $\oX$.
Denote by $j_i: Y_i \rightarrow \oX$ the inclusions and let 
$$
\cN_{\oX,Y_1,\ldots,Y_n} \subseteq 
\cN_{\oX}\oplus\bigoplus_{i=1}^n j_{i*}\cN_{Y_i}
$$
be the subsheaf consisting of those tuples of local sections
$(s,t_1,\ldots,t_n)$ such that $s\vert_{Y_i}=p_i(t_i)$ for all $i$.  
Then $\cN_{\oX,Y_1,\ldots,Y_n}$ is a sheaf on $\oX$, 
endowed with a morphism
$$
\eta:\cT_{\bP(V)}\vert_{\oX} \rightarrow 
\cN_{\oX,Y_1,\ldots,Y_n}
$$
that sends any tangent vector to the collection of the
corresponding normal vectors to $\oX,Y_1,\ldots,Y_n$.
The first projection 
$\cN_{\oX,Y_1,\ldots,Y_n}\rightarrow \cN_{\oX}$
is surjective, with kernel the direct sum of the $\cN_{Y_i/\oX}$.
In addition, $\eta$ fits into a commutative diagram
$$
\CD
0 @>>> \cT_{\oX} @>>> \cT_{\bP(V)}\vert_{\oX} @>>>
\cN_{\oX}  @>>> 0\\
& & @VVV @V{\eta}VV @V{\id}VV \\
0 @>>> \bigoplus_{i=1}^n \cN_{Y_i/\oX} @>>> 
\cN_{\oX,Y_1,\ldots,Y_n} @>>> \cN_{\oX}  @>>> 0.\\
\endCD
$$
Together with the exact sequence
$$
0 \rightarrow \cT_{\oX}(-\log \partial\oX) \rightarrow
\cT_{\oX} \rightarrow \bigoplus_{i=1}^n \cN_{Y_i/\oX} 
\rightarrow 0,
$$
this yields an exact sequence
$$
\CD
0 @>>> \cT_{\oX}(-\log\partial\oX)
@>>> \cT_{\bP(V)}\vert_{\oX} @>{\eta}>>
\cN_{\oX,Y_1,\ldots,Y_n} @>>> 0.
\endCD
$$
Since $H^1(\oX,\cT_{\oX}(-\log\partial\oX))=0$, we obtain 
a surjection
$$
H^0(\eta): H^0(\oX,\cT_{\bP(V)}) \rightarrow
H^0(\oX,\cN_{\oX,Y_1,\ldots,Y_n}) 
= T_{(\oX,Y_1,\ldots,Y_n)} \Hilb.
$$
Note that $f$ factors as the restriction
$$
\rho:\End(V)\rightarrow H^0(\oX,\cT_{\bP(V)}),
$$
followed by $H^0(\eta)$. So it remains to check surjectivity 
of $\rho$. 

For this, we use the exact sequence
$$
0 \rightarrow \cO_{\oX} \rightarrow V\otimes_k \cO_{\oX}(1)
\rightarrow \cT_{\bP(V)}\vert_{\oX} \rightarrow 0,
$$
and the vanishing of $H^1(\oX,\cO_{\oX})$ (since $\oX$ is
nonsingular, projective and rational). This yields a 
surjection 
$$
V\otimes H^0(\oX,\cO_{\oX}(1)) \rightarrow 
H^0(\oX, \cT_{\bP(V)}).
$$
Since the restriction $V^*\rightarrow H^0(\oX,\cO_{\oX}(1))$
is surjective as well, the composition 
$\End(V) = V\otimes V^*\rightarrow
H^0(\oX, \cT_{\bP(V)}\vert_{\cO_{\oX}})$
is surjective as desired.
\end{proof}

Consider the action of $G$ on the Grassmannian $\Grass(\fg)$
of subspaces of its Lie algebra. If $H$ is a spherical subgroup
of $G$ with Lie algebra $\fh$, then the isotropy group of 
$\fh\in \Grass(\fg)$ is the normalizer $N_G(\fh)$. The latter 
equals $N_G(H)$, by \cite[\S 5]{Brion_Pauer87}. 
So Theorem \ref{fin} implies 

\begin{corollary}\label{lun}
There exist only finitely many conjugacy classes of spherical
subgroups $H\subseteq G$ having finite index in their 
normalizer.
\end{corollary}

(Note that there exist infinitely many conjugacy classes of
spherical subgroups of $G$, for example, those containing $U$ 
as a subgroup of finite index.) For $G$ of type $A$, Corollary 
\ref{lun} is a consequence of Luna's classification of spherical
varieties by combinatorial invariants \cite{Luna01}. Extending
this classification to arbitrary $G$ would give another proof
of that corollary.

Another direct consequence of Theorem \ref{fin} is 

\begin{corollary}\label{fin1}
Any finite-dimensional $G$-module $V$ contains only finitely many 
multiplicity-free subvarieties, up to the action of $\GL(V)^G$.
\end{corollary}

\begin{proof}
Let $X$ be a multiplicity-free subvariety of $V$, with open orbit
$G\cdot x$. By Theorem \ref{fin}, we may fix the isotropy group 
$G_x = H$.
Write
$$
V=\bigoplus_{\lambda\in F} M_{\lambda}\otimes_k V(\lambda)
$$
where $F$ is a finite subset of $\Lambda^+$, and 
$M_{\lambda}=\Hom^G(V(\lambda),V)$. Then
$$
\GL(V)^G\simeq \prod_{\lambda\in F} \GL(M_{\lambda}),
$$
and 
$x\in V^H=\bigoplus_{\lambda\in F} M_{\lambda}\otimes V(\lambda)^H$
where every $V(\lambda)^H$ is at most one-dimensional. Thus,
$$
x=\sum_{\lambda\in F} m_{\lambda}\otimes x_{\lambda}
$$
where every nonzero $x_{\lambda}$ is a generator of $V(\lambda)^H$. 
It follows that there are only finitely many $x$ up to the action 
of $\GL(V)^G$. So the same holds for $X=\overline{G\cdot x}$.
\end{proof}

Together with Theorem \ref{rep3}, this implies

\begin{corollary}\label{fin2}
Given a submonoid $\cS\subseteq\Lambda^+$, there exist only finitely 
many isomorphism classes of multiplicity-free $G$-varieties with 
weight monoid $\cS$. Equivalently, the moduli scheme $\M_{\cS}$ 
contains only finitely many $T_{\ad}$-orbits.
\end{corollary}


\end{document}